\documentclass[11pt, english, twoside]{article}
\usepackage[margin=3cm]{geometry}
\usepackage{amsbsy,amsmath,amsthm,amssymb}

\usepackage{fourier}
\DeclareMathAlphabet{\mathcal}{OMS}{cmsy}{m}{n}

\usepackage[dvipsnames]{xcolor}
\usepackage{bm}
\usepackage{esint}
\usepackage{babel}
\usepackage[colorlinks=true, citecolor = blue]{hyperref}

\usepackage{indentfirst}
\usepackage{changepage, cases}	
\usepackage{setspace}
\usepackage{indentfirst}
\usepackage{graphicx}
\usepackage{calc}
\usepackage{orcidlink}
\usepackage{booktabs,multirow}
\usepackage{makecell}

\usepackage{fancyhdr}
\allowdisplaybreaks

\theoremstyle{plain}
\newtheorem{theorem}{Theorem}[section]

\newtheorem{proposition}[theorem]{Proposition}

\theoremstyle{remark}
\newtheorem{remark}[theorem]{Remark}

\theoremstyle{definition}

\newenvironment{proof of theorem 1.1}{{\noindent \em Proof of Theorem 1.1.}}{\hfill $\Box$\par}
\newenvironment{proof of theorem 1.2}{{\noindent \em Proof of Theorem 1.2.}}{\hfill $\Box$\par}

\DeclareSymbolFont{EulerExtension}{U}{euex}{m}{n}
\DeclareMathSymbol{\euintop}{\mathop} {EulerExtension}{"52}
\DeclareMathSymbol{\euointop}{\mathop} {EulerExtension}{"48}

\begin{document}
	\title{On the multivariate hypergeometric function $F_K^{(L)}$, I}
	\author{
		Min-Jie Luo \orcidlink{0000-0001-7433-4490} $^{\rm 1}$\thanks{Corresponding author}~, 
		Shu-Ying Shen $^{\rm 2}$} 
	\date{}
	\maketitle
	\begin{center}\small
		$^{1}$\emph{School of Mathematics and Statistics, Donghua University,\\ 
			Shanghai 201620, People's Republic of China.}\\
			E-mail: \url{mathwinnie@live.com}, \url{mathwinnie@dhu.edu.com}
	\end{center}
	\begin{center}\small
		$^{2}$\emph{School of Mathematics and Statistics, Donghua University,\\ 
			Shanghai 201620, People's Republic of China.}\\
		E-mail: \url{15207010215@163.com}
	\end{center}
	
	
\begin{abstract}
	In recent years, a new multivariate hypergeometric series, known as the comb function, has attracted attention in physics, particularly in the study of conformal blocks in conformal field theory. This function can be regarded as a further generalization of Appell's function $F_2$ and Saran's function $F_K$. However, except for a few known formulas, little is known about its properties. In this paper, we study its fundamental properties from a mathematical perspective, such as the region of convergence, reduction formulas, an Euler-type integral representation, and decomposition formulas. Some directions for future research and connections with other multivariate hypergeometric series originating from physics are also mentioned. \\
		
		
\end{abstract}



\section{Introduction}\label{Introduction}

In 1888, Lauricella \cite{Lauricella-1893} introduced 14 hypergeometric series in three variables, namely $F_1,\cdots,F_{14}$. Among them, $F_1$, $F_2$, $F_5$ and $F_9$ are what we now know as the Lauricella functions $F_A$, $F_B$, $F_C$, $F_D$, while the remaining $10$ fell into oblivion for a long time. Beginning in 1954, Saran published a series of papers (\cite{Saran-1954}, \cite{Saran-1955a}, \cite{Saran-1955b} and \cite{Saran-1956}), which brought these $10$ functions back into view. Therefore, in the literature, these functions are also called \emph{Lauricella-Saran functions}.

Although they all satisfy ``hypergeometricity,'' subtle differences in the structure of their coefficients often lead to vastly different properties of the functions, and therefore each function needs to be studied individually. Around 2021, Luo (the first author of this paper) and Raina \cite{Luo-Raina-2021} realized that Saran's $F_K$ function possesses a certain ``potential'' both in theory and in applications. This function can be defined as (see \cite[p. 2, Eq. (2)]{Luo-Raina-2021}; see also \cite[p. 42, Eq. (4)]{Srivastava-Karlsson-Book-1985} and \cite[p. 67, Eq. (29) ]{Srivastava-Manocha-Book-1984}):
\begin{align}\label{Def-SaranFK}
	F_K\big[\alpha_1, \alpha_2, \alpha_2, \beta_1, \beta_2, \beta_1;\gamma_1, \gamma_2, \gamma_3; x, y, z\big]
	:=\sum_{m,n,p=0}^{\infty} 
	\frac{(\alpha_1)_m (\alpha_2)_{n+p} (\beta_1)_{m+p} (\beta_2)_n}{(\gamma_1)_m (\gamma_2)_n (\gamma_3)_p} 
	\frac{x^m}{m!} \frac{y^n}{n!} \frac{z^p}{p!}, 
\end{align}
where $(x,y,z)\in \mathbb{D}:=\left\{(x,y,z)\in\mathbb{C}^3;|x|<1,|y|<1,|z|<(1-|x|)(1-|y|)\right\}$. Note that $\mathbb{D}$ is a complete Reinhardt domain \cite[p. 1]{Luo et al-2022}. Through simple calculations, we can transform \eqref{Def-SaranFK} into the following equivalent forms \cite{Hang-Luo-2025}:
\begin{align}
	&F_K\big[\alpha_1,\alpha_2,\alpha_2,\beta_1,\beta_2,\beta_1;\gamma_1,\gamma_2,\gamma_3;x,y,z\big]\notag\\
	&\hspace{1cm}=\sum_{p=0}^{\infty}\frac{(\alpha_2)_p(\beta_1)_p}{p!(\gamma_3)_p}
	{}_{2}F_{1}\left[\begin{matrix}
		\beta_1+p,\alpha_1\\
		\gamma_1
	\end{matrix}; x\right]
	{}_{2}F_{1}\left[\begin{matrix}
		\alpha_2+p,\beta_2\\
		\gamma_2
	\end{matrix}; y\right]z^p\label{FK-Gauss2F1}\\
	&\hspace{1cm}=\sum_{m=0}^{\infty}\frac{(\alpha_1)_m(\beta_1)_m}{m!(\gamma_1)_m}F_2\big[\alpha_2,\beta_2,\beta_1+m;\gamma_2,\gamma_3;y,z\big]x^m,\label{FK-AppellF2}
\end{align}
where
\begin{equation}\label{Def-Gauss2F1}
{}_{2}F_{1}\left[\begin{matrix}
	a,b\\
	c
	\end{matrix}; z\right]:=\sum_{n=0}^{\infty}\frac{(a)_n(b)_n}{(c)_n}\frac{z^n}{n!}\quad (|z|<1)
\end{equation}
is the Gauss hypergeometric series \cite[p. 383, Chapter 15]{Olver et al-NISThandbook-2010} and
\[
F_2\big[a,b,b';c,c';x,y\big]:=\sum_{m,n=0}^{\infty}\frac{(a)_{m+n}(b)_m(b')_n}{(c)_m(c')_n}\frac{x^m}{m!}\frac{y^n}{n!}\quad (|x|+|y|<1)
\]
is the second Appell hypergeometric series \cite[p. 23]{Srivastava-Karlsson-Book-1985}. 

In their work \cite{Luo-Raina-2021}, a (single) Mellin--Barnes contour integral representation of $F_K$ was established, from which its analytic continuation and large-variable asymptotic behaviors were obtained. In addition, by means of fractional integration by parts, they also obtained an Erd\'{e}lyi-type integral. The study of integrals of this type originated in Erd\'{e}lyi's paper \cite{Erdelyi-1939b}, while the fractional integration by parts technique grew out of a series of studies on a problem in Riemann's 1857 paper \cite{Riemann-1857}.

Paper \cite{Luo-Raina-2021} gave rise to two lines of research. First, in the context of Erd\'{e}lyi-type integrals, Luo \emph{et al.} \cite{Luo et al-2022} extended the ideas in \cite{Luo-Raina-2021} and obtained a more elegant integral. Then, this work was further extended by Guo and Luo \cite{Guo-Luo-2026}. Second, regarding the asymptotics of functions, with the help of the theory of asymptotic expansions for Mellin integrals, Hang and Luo \cite{Hang-Luo-2025} characterized the more complicated asymptotic behavior of the $F_K$ function, namely the case when two of its variables become simultaneously large. In addition, as a by-product, Hang and Luo established many asymptotic properties of Humbert's function $\Psi_1$ in \cite{Hang-Luo-2025}, which in turn motivated several other studies on Humbert functions (see \cite{Hang-Luo-2024}, \cite{Hang et al-2025}, \cite{Hang et al-2026a} and \cite{Hang et al-2026b}). The above developments reflect the diversity and breadth of the theory of a single multivariate hypergeometric function.

There are also some other interesting results concerning the $F_K$ function that have been obtained in recent years. For example, Antonova \emph{et al.} \cite{Antonova-Dmytryshyn-Goran-2023} studied ratios of the $F_K$ function via the theory of branched continued fractions. For further developments in this line of research, see \cite{Dmytryshyn-Goran-2024} and \cite{Dmytryshyn-Goran-2026}.

In this paper, we study an important class of $L$-dimensional generalizations of Saran's $F_K$ function \eqref{Def-SaranFK}. These generalizations are defined by the following $L$-fold power series \cite[p. 3, Eq. (8)]{Luo-Raina-2021}:
\begin{align}\label{Def-MultivariateFK}
	&F_K^{(L)}\left[\begin{matrix}
		a_1, b_1, \cdots, b_{L-1}, a_2 \\ 
		c_1, \cdots, c_L
		\end{matrix}; z_1,\cdots,z_L\right]\notag\\
	&\hspace{1cm}:=\sum_{n_1, \cdots, n_L=0}^{\infty} \frac{(a_1)_{n_1}(b_1)_{n_1+n_2}\cdots(b_{L-1})_{n_{L-1}+n_L}(a_2)_{n_L}}{(c_1)_{n_1}\cdots(c_L)_{n_L}} \frac{z_1^{n_1}}{n_1!} \cdots \frac{z_L^{n_L}}{n_L!}.
\end{align}
The region of convergence of series \eqref{Def-MultivariateFK} will be discussed in detail in Section \ref{Sect-2}. Obviously, we have
\begin{align*}
F_K^{(1)}\left[\begin{matrix}
	a_1, a_2 \\ 
	c_1
\end{matrix}; z_1\right]&={}_{2}F_{1}\left[\begin{matrix}
a_1, a_2 \\ 
c_1
\end{matrix}; z_1\right],\notag\\
F_K^{(2)}\left[\begin{matrix}
	a_1, b_1, a_2 \\ 
	c_1, c_2
\end{matrix}; z_1,z_2\right]
&=F_2\big[b_1,a_1,a_2;c_1,c_2;z_1,z_2\big],\notag\\
F_K^{(3)}\left[\begin{matrix}
	a_1, b_1, b_2, a_2 \\ 
	c_1, c_2, c_3
\end{matrix}; z_1,z_2,z_3\right]
&=F_K\big[a_1, b_2, b_2, b_1, a_2, b_1;c_1, c_3, c_2; z_1, z_3, z_2\big].
\end{align*}

The $F_K^{(L)}$ function was first introduced by Belitsky in his 2018 paper \cite{Belitsky-2018}, where it was used to provide an explicit representation of the multichannel conformal blocks obtained by resumming all small fermion\textcolor{orange}{--}antifermion pairs (conformal descendants) in the flux-tube description of $N$-sided NMHV scattering amplitudes. Later, Rosenhaus \cite{Rosenhaus-2019} used this function while studying the $n$-point conformal blocks in the comb channel and referred to it as the \emph{comb function}. In the appendix of \cite{Rosenhaus-2019}, he discussed some basic properties of this function, such as splitting identities, the differential equations it satisfies, and the Euler-type integral representation. Ferrando \emph{et al.} \cite{Ferrando et al-2026} used a \emph{rescaled version} of this function in evaluating the conformal integral associated with the $n$-point conformal partial wave in the comb channel. For its other physical applications, interested readers may refer to \cite{Ammon et al-2025}, \cite{Bhattacharyya et al-2025}, \cite{Pavlov-2023}, etc.

Throughout the paper, boldface letters always denote vectors of dimension $L$; for instance, $\mathbf{n}:=(n_1,\cdots,n_L)\in\mathbb{Z}_{\geq0}^{L}$ and $\mathbf{z}:=(z_1,\cdots,z_L)\in\mathbb{C}^L$. Let $\mathbf{e}_i$ denote the $i$-th standard unit vector in $\mathbb{R}^{L}$, and let $\langle\cdot,\cdot\rangle$ denote the inner product on $\mathbb{R}^{L}$. We write $\mathbf{z}^{\mathbf{n}}:=z_1^{n_1}\cdots z_L^{n_L}$ and $\mathbf{n}!:=n_1!\cdots n_L!$. In addition,  
\[
\sum_{\mathbf{n}=\mathbf{0}}^{\infty}\Omega(\mathbf{n})
\quad\text{means}\quad
\sum_{n_1,\cdots,n_L=0}^{\infty}\Omega(n_1,\cdots,n_L). 
\]
With the above notation, we often write \eqref{Def-MultivariateFK} as
\[
F_K\left[\begin{matrix}
	a_1, b_1, \cdots, b_{L-1}, a_2 \\ 
	c_1, \cdots, c_L
\end{matrix}; \mathbf{z}\right]
=\sum_{\mathbf{n}=\mathbf{0}}^{\infty} \frac{(a_1)_{n_1}(b_1)_{n_1+n_2}\cdots(b_{L-1})_{n_{L-1}+n_L}(a_2)_{n_L}}{(c_1)_{n_1}\cdots(c_L)_{n_L}} 
\frac{\mathbf{z}^{\mathbf{n}}}{\mathbf{n}!}.
\]

The paper is organized as follows. In Section \ref{Sect-2}, we derive some reduction formulas for $F_{K}^{(L)}$ and use them to characterize its region of convergence. In Section \ref{Sect-3}, we use Rosenhaus's double integral representation to derive a multiple integral representation for $F_K^{(L)}$ whose integrand is an elementary function. Subsequently, in Section \ref{Sect-4}, we apply symbolic operators to obtain some decomposition formulas for $F_K^{(L)}$. Finally, in Section \ref{Sect-5}, we outline some directions for future research and discuss connections between $F_K^{(L)}$ and another recently introduced multivariable hypergeometric function, namely the \emph{necklace function}.

\section{Reduction formulas}\label{Sect-2}

Reduction formulas play an important role in the study of hypergeometric functions. They help us identify the parameter combinations under which a hypergeometric function can be represented by simpler special functions or even elementary functions. In addition to simplifying calculations, reduction formulas are powerful tools for determining the convergence regions of many hypergeometric functions (see \cite[p. 118, Section 4.6]{Srivastava-Karlsson-Book-1985} and \cite[p. 212, Section 2]{Karlsson-1986}). For example, if we let $\gamma_1=\alpha_1$ and $\gamma_2=\beta_2$ in \eqref{FK-Gauss2F1}, then the two ${}_{2}F_{1}$ functions involved can be evaluated via the reduction formula
\begin{equation}\label{ReductionF-2F1}
{}_{2}F_{1}\left[\begin{matrix}
	a,b\\
	b
\end{matrix}; z\right]=(1-z)^{-a}\quad (|z|<1)
\end{equation}
and thus we obtain the following reduction formula:
\begin{align*}
	&F_K\big[\alpha_1,\alpha_2,\alpha_2, \beta_1, \beta_2, \beta_1; \alpha_1, \beta_2, \gamma_3; x,y,z\big]\notag\\
	&\hspace{1cm}=(1-x)^{-\beta_1}(1-y)^{-\alpha_2}
	{}_{2}F_{1}\left[\begin{matrix}
		\alpha_2, \beta_1\\
		\gamma_3
	\end{matrix};\frac{z}{(1-x)(1-y)}\right],
\end{align*}
where $|x|<1$, $|y|<1$, and $|z|<(1-|x|)(1-|y|)$. Clearly, the convergence condition here is obtained by applying \eqref{ReductionF-2F1} twice and the series definition \eqref{Def-Gauss2F1}, and it agrees with the convergence condition for $F_K$ under the general definition \eqref{Def-SaranFK}. However, this convergence condition is not easy to see from the general definition \eqref{Def-SaranFK} or its equivalent representations \eqref{FK-Gauss2F1} and \eqref{FK-AppellF2}.

Next, let us derive some reduction formulas for $F_K^{(L)}$. 
\begin{theorem}\label{Prop-ReductionF}
	Let $L\geq 2$. Also, let $\Phi_L(z_L):=1-z_L$ and 
	\[
	\Phi_k(z_k,\cdots,z_L):=1-\dfrac{z_k}{\Phi_{k+1}(z_{k+1},\cdots,z_L)}
	\quad (k=1,\cdots,L-1).
	\] 
	For simplicity, whenever no confusion can arise, we write $\Phi_k\equiv \Phi_k(z_k,\cdots,z_L)$. Then we have
	\begin{equation}\label{Prop-ReductionF-1}
		F_K^{(L)}\left[\begin{matrix}
			a_1,b_1,b_2,\cdots,b_{L-1},a_2\\
			c_1,b_2,\cdots,b_{L-1},a_2
		\end{matrix};\mathbf{z}
		\right]
		=\Phi_2^{-b_1}\cdots\Phi_L^{-b_{L-1}}~
		{}_{2}F_{1}\left[\begin{matrix}
			a_1,b_1\\
			c_1
		\end{matrix};1-\Phi_1\right],
	\end{equation}
	where $|1-\Phi_k|<1$ $(k=1,\cdots,L)$. 
	
	In particular, if we further set $c_1=b_1$ in \eqref{Prop-ReductionF-1}, then we obtain  
	\[
		F_K^{(L)}\left[\begin{matrix}
			a_1,b_1,\cdots,b_{L-1},a_2\\
			b_1,\cdots,b_{L-1},a_2
		\end{matrix};\mathbf{z}
		\right]
		=\Phi_1^{-a_1}\Phi_2^{-b_1}\cdots\Phi_L^{-b_{L-1}},
	\]
	where $|1-\Phi_k|<1$ $(k=1,\cdots,L)$.
\end{theorem}
\begin{proof}
	Using $(a)_{m+n}=(a)_n(a+n)_m$ and the series definition \eqref{Def-MultivariateFK}, we have
	\begin{align}\label{Prop-ReductionF-Proof-1}
	F_K^{(L)}\left[\begin{matrix}
		a_1,b_1,b_2,\cdots,b_{L-1},a_2\\
		c_1,b_2,\cdots,b_{L-1},a_2
	\end{matrix};\mathbf{z}\right]
	&=\sum_{\mathbf{n}=\mathbf{0}}^\infty
	\frac{(a_1)_{n_1}(b_1)_{n_1}}{(c_1)_{n_1}}
	(b_1+n_1)_{n_2}\cdots (b_{L-1}+n_{L-1})_{n_L}
	\frac{\mathbf{z}^{\mathbf{n}}}{\mathbf{n}!}\notag\\
	&=\sum_{\mathbf{n}=\mathbf{0}}^\infty
	\frac{(a_1)_{n_1}(b_1)_{n_1}}{(c_1)_{n_1}}
	\frac{z_1^{n_1}}{n_1!}
	\cdot
	\prod_{j=2}^{L}
	(b_{j-1}+n_{j-1})_{n_j}\frac{z_j^{n_j}}{n_j!}.
	\end{align}
	So the result can be obtained by preforming summation over the indices $n_L$, $n_{L-1}$, $\cdots$, $n_1$ sequentially. 
	
	First, summing over the index $n_L$ and using the binomial theorem gives
	\begin{equation}\label{Prop-ReductionF-Proof-2}
	\sum_{n_L=0}^{\infty}(b_{L-1}+n_{L-1})_{n_L}\frac{z_L^{n_L}}{n_L!}
	=\Phi_L^{-b_{L-1}-n_{L-1}}, 
	\end{equation}
	where $\Phi_L=1-z_L$ satisfies $|1-\Phi_L|=|z_L|<1$. Substituting \eqref{Prop-ReductionF-Proof-2} into \eqref{Prop-ReductionF-Proof-1}, we obtain
	\begin{align}\label{Prop-ReductionF-Proof-3}
		&F_K^{(L)}\left[\begin{matrix}
			a_1,b_1,b_2,\cdots,b_{L-1},a_2\\
			c_1,b_2,\cdots,b_{L-1},a_2
		\end{matrix};\mathbf{z}\right]\notag\\
		&\hspace{0.5cm}=\Phi_L^{-b_{L-1}}\sum_{n_1,\dots,n_{L-1}=0}^\infty
		\frac{(a_1)_{n_1}(b_1)_{n_1}}{(c_1)_{n_1}}\frac{z_1^{n_1}}{n_1!}
		\Phi_{L}^{-n_{L-1}}
		\prod_{j=2}^{L-1}
		(b_{j-1}+n_{j-1})_{n_j}\frac{z_j^{n_j}}{n_j!},
	\end{align}
	where $|1-\Phi_L|<1$. Next, summing over the index $n_{L-1}$ gives 
	\begin{equation}\label{Prop-ReductionF-Proof-4}
		\sum_{n_{L-1}=0}^{\infty}\Phi_{L}^{-n_{L-1}}
		(b_{L-2}+n_{L-2})_{n_{L-1}}\frac{z_{L-1}^{n_{L-1}}}{n_{L-1}!}
		=\left(1-\frac{z_{L-1}}{\Phi_L}\right)^{-b_{L-2}-n_{L-2}}
		=\Phi_{L-1}^{-b_{L-2}-n_{L-2}},
	\end{equation}
	where $\Phi_{L-1}=1-\dfrac{z_{L-1}}{\Phi_L}$ satisfies $|1-\Phi_{L-1}|<1$. Substituting \eqref{Prop-ReductionF-Proof-4} into \eqref{Prop-ReductionF-Proof-3} we obtain
	\begin{align*}
		&F_K^{(L)}\left[\begin{matrix}
			a_1,b_1,b_2,\cdots,b_{L-1},a_2\\
			c_1,b_2,\cdots,b_{L-1},a_2
		\end{matrix};\mathbf{z}\right]\\
		&\hspace{0.5cm}=\Phi_L^{-b_{L-1}}\Phi_{L-1}^{-b_{L-2}}\sum_{n_1,\dots,n_{L-2}=0}^\infty
		\frac{(a_1)_{n_1}(b_1)_{n_1}}{(c_1)_{n_1}}\frac{z_1^{n_1}}{n_1!}
		\Phi_{L-1}^{-n_{L-2}}\prod_{j=2}^{L-2}
		(b_{j-1}+n_{j-1})_{n_j}\frac{z_j^{n_j}}{n_j!},
	\end{align*}
	where $|1-\Phi_L|<1$ and $|1-\Phi_{L-1}|<1$. By continuing this process, we finally obtain
	\begin{align}\label{Prop-ReductionF-Proof-5}
		F_K\left[\begin{matrix}
			a_1,b_1,b_2,\cdots,b_{L-1},a_2\\
			c_1,b_2,\cdots,b_{L-1},a_2
		\end{matrix};\mathbf{z}\right]
	&=\Phi_L^{-b_{L-1}}\cdots \Phi_3^{-b_2}
		\sum_{n_1,n_2=0}^\infty
		\frac{(a_1)_{n_1}(b_1)_{n_1}}{(c_1)_{n_1}}\frac{z_1^{n_1}}{n_1!}
		\Phi_{3}^{-n_2}(b_1+n_1)_{n_2}
		\frac{z_2^{n_2}}{n_2!}\notag\\
	&=\Phi_L^{-b_{L-1}}\cdots \Phi_3^{-b_2}\sum_{n_1=0}^\infty
		\frac{(a_1)_{n_1}(b_1)_{n_1}}{(c_1)_{n_1}}
		\frac{z_1^{n_1}}{n_1!}
		\left(1-\frac{z_2}{\Phi_3}\right)^{-b_1-n_1}\notag\\
	&=\Phi_L^{-b_{L-1}}\cdots \Phi_2^{-b_1}\sum_{n_1=0}^\infty
		\frac{(a_1)_{n_1}(b_1)_{n_1}}{(c_1)_{n_1}}
		\frac{(z_1/\Phi_2)^{n_1}}{n_1!},
	\end{align}
	where $|1-\Phi_k|<1$ $(k=1,\cdots,L)$. Now the result \eqref{Prop-ReductionF-1} follows by interpreting the last series in \eqref{Prop-ReductionF-Proof-5} as the Gauss hypergeometric function and noting that $1-\Phi_1=\dfrac{z_1}{\Phi_2}$.
\end{proof}
\begin{remark}
Letting $L=2$ in \eqref{Prop-ReductionF-1}, it becomes \cite[p. 306, Eq. (109)]{Srivastava-Karlsson-Book-1985}
	\[
	F_2\big[b_1,a_1,a_2;c_1,a_2;z_1,z_2\big]
	=(1-z_2)^{-b_1}{}_{2}F_{1}\left[\begin{matrix}
		a_1,b_1\\
		c_1
	\end{matrix};\frac{z_1}{1-z_2}\right].
	\]
	When $L=3$ in \eqref{Prop-ReductionF-1}, we obtain
	\[
	F_K^{(3)}\left[\begin{matrix}
		a_1,b_1,b_2,a_2\\
		c_1,b_2,a_2
	\end{matrix};z_1,z_2,z_3
	\right]
	=(1-z_2-z_3)^{-b_1}(1-z_3)^{b_1-b_2} 
	{}_{2}F_{1}\left[\begin{matrix}
		a_1,b_1\\
		c_1
	\end{matrix};\frac{z_1(1-z_3)}{1-z_2-z_3}\right].
	\]
\end{remark}

By Theorem \ref{Prop-ReductionF}, we can give an exact characterization of the convergence region of $F_K^{(L)}$. 

\begin{proposition}\label{Prop-ReductionF-ConvergenceC}
	Let $r_k=|z_k|\geq0$ $(k=1,\cdots,L)$. Then the multiple series defined by \eqref{Def-MultivariateFK} is absolutely convergent in 
	\begin{equation}\label{Prop-ReductionF-ConvergenceC-1}
	\mathbb{D}_L:=\left\{\mathbf{z}\in\mathbb{C}^L;~ \Phi_k(r_k,\cdots,r_L)>0,~k=1,\cdots,L\right\},
	\end{equation}
	where $\Phi_k$ is defined in Theorem \ref{Prop-ReductionF}. 
\end{proposition}
\begin{remark}
~
\begin{itemize}
	\item[(1)] In 2020, the first author of the present paper supervised the undergraduate student Jing-Feng Quan at Donghua University to investigate the region of convergence of series \eqref{Def-MultivariateFK} via reduction formulas \cite{Quan-2020}. At that time, they considered the case where all parameters are equal to $1$, and the resulting region of convergence agrees with Proposition \ref{Prop-ReductionF-ConvergenceC}.
	\item[(2)] Regarding the determination of the convergence domains of multivariable hypergeometric series, another well-known method is to use Horn's results \cite{Horn-1889}. This method is concise and effective when dealing with double hypergeometric series. However, as Srivastava and Karlsson \cite[p. 268]{Srivastava-Karlsson-Book-1985} mention, although Horn's method can be generalized to $n$ variables, the expressions involved quickly become extremely complicated as $n$ grows. This demonstrates the advantage of our method for determining the convergence region via reduction formulas. In Appendix \ref{Appendix-1}, we use Horn's method to obtain the region of convergence of $F_K^{(4)}$ and show that it is the same as the result derived from Proposition \ref{Prop-ReductionF-ConvergenceC}.
\end{itemize}
\end{remark}

Obviously, for $L=1$, the region \eqref{Prop-ReductionF-ConvergenceC-1} becomes 
\[
\mathbb{D}_1=\big\{z_1\in\mathbb{C};~r_1<1\big\}
=\big\{z_1\in\mathbb{C};~|z_1|<1\big\},
\] 
which is the region of convergence of the Gauss hypergeometric series ${}_{2}F_{1}$. For $L=2$, the region \eqref{Prop-ReductionF-ConvergenceC-1} becomes 
\[
\mathbb{D}_2=\left\{(z_1,z_2)\in\mathbb{C}^2;~r_2<1,~\frac{r_1}{1-r_2}<1\right\}
=\big\{(z_1,z_2)\in\mathbb{C}^2;~|z_1|+|z_2|<1\big\},
\]
which is the standard region of convergence of the Appell function $F_2$. For $L=3$, we have 
\begin{align*}
\mathbb{D}_3&=\left\{(z_1,z_2,z_3)\in\mathbb{C}^3;~ r_3<1,~ \frac{r_2}{1-r_3}<1,~\frac{r_1}{1-\frac{r_2}{1-r_3}}<1\right\}\\
&=\big\{(z_1,z_2,z_3)\in\mathbb{C}^3;~ |z_1|<1,~ |z_3|<1, ~ |z_2|<(1-|z_1|)(1-|z_3|)\big\}.
\end{align*}
This is the convergence region of Saran's $F_K$ function. For $L=4$, we note that the conditions $\Phi_k(r_k,\cdots,r_4)>0, k=1,2,3,4$, are equivalent to
\[
	r_4<1,~r_1<1,~\frac{r_3}{1-r_4}+\frac{r_2}{1-r_1}<1.
\]
Thus 
\begin{equation}\label{Cond-FK4}
\mathbb{D}_4=\left\{(z_1,z_2,z_3,z_4)\in\mathbb{C}^4;~ |z_1|<1,~ |z_4|<1,~ \frac{|z_2|}{1-|z_1|}+\frac{|z_3|}{1-|z_4|}<1 \right\}, 
\end{equation}
which agrees with the conditions stated in \cite[p. 11]{Guo-Luo-2026}.

We have mentioned in Section \ref{Introduction} that $\mathbb{D}_L$ is a complete Reinhardt domain. A set $S\subseteq\mathbb{C}^L$ is called a complete Reinhardt domain (or a complete circular domain) if $\mathbf{z}\in S$ implies that $(\mu_1 z_1,\cdots,\mu_L z_L)\in S$ for all $\mu_j\in \mathbb{C}$ with $|\mu_j|\leq 1$, $j=1,\cdots,L$ (see \cite[pp. 104--105, Definition 2.3.12]{Krantz-Book-2001}). The set is said to be logarithmically convex if 
\[
\log\|S\|:=\big\{(\log|z_1|,\cdots,\log|z_L|);~ \mathbf{z}=(z_1,\cdots,z_L)\in S\big\}\subseteq \mathbb{R}^L
\]
is convex (see \cite[p. 105, Definition 2.3.14]{Krantz-Book-2001}).
\begin{proposition}\label{Prop-ReductionF-Rienhardt}
	The set~~$\mathbb{D}_L$ defined by \eqref{Prop-ReductionF-ConvergenceC-1} is a complete Reinhardt domain and is logarithmically convex. 
\end{proposition}
\begin{proof}
	The result follows from Proposition 2.3.15 in \cite[p. 205]{Krantz-Book-2001}. However, given the complexity of $\mathbb{D}_L$, it may be helpful to readers unfamiliar with several complex variables to provide a direct verification that does not rely on this proposition.

	It is easy to verify by induction that 
	\begin{equation}\label{Prop-ReductionF-Rienhardt-Proof-1}
	0\leq s_j\leq r_j~(j=1,\cdots,L)
	~\implies~
	\Phi_k(s_k,\cdots,s_L)\geq \Phi_k(r_k,\cdots,r_L)>0 
	~(k=1,\cdots,L).
	\end{equation}
	When $k=L$, we have
	$\Phi_L(s_L)\geq\Phi_L(r_L)>0$. Assume  $\Phi_{k+1}(s_{k+1},\cdots,s_L)\geq \Phi_{k+1}(r_{k+1},\cdots,r_L)>0$. Since $s_k\leq r_k$, we have
	\[
	\Phi_k(s_k,\cdots,s_L)
	=1-\frac{s_k}{\Phi_{k+1}(s_{k+1},\cdots,s_L)}
	\geq1-\frac{r_k}{\Phi_{k+1}(r_{k+1},\cdots,r_L)}
	=\Phi_{k}(r_k,\cdots,r_L)>0.
	\]
	Let $\mathbf{z}\in \mathbb{D}_L$, and let $\bm{\mu}\in\mathbb{C}^L$ such that $|\mu_j|\leq 1$ $(j=1,\cdots,L)$. According to \eqref{Prop-ReductionF-Rienhardt-Proof-1}, we have 
	\[
	\Phi_k(|\mu_k z_k|,\cdots,|\mu_L z_L|)\geq \Phi_k(|z_k|,\cdots,|z_L|)>0 
	~(k=1,\cdots,L),
	\]
	which implies $(\mu_1 z_1,\cdots,\mu_L z_L)\in\mathbb{D}_L$. Therefore $\mathbb{D}_L$ is complete Reinhardt. 
	
	To prove that $\mathbb{D}_L$ is logarithmically convex, we show that
	\[
	\log \|\mathbb{D}_L\|=\big\{\mathbf{x}\in\mathbb{R}^L;~ \Phi_k(\mathrm{e}^{x_k},\cdots,\mathrm{e}^{x_L})>0,~k=1,\cdots,L\big\}\subset\mathbb{R}^L
	\] 
	is convex, which is equivalent to prove that for any $\mathbf{x},\mathbf{y}\in \log \|\mathbb{D}_L\|$ and any $\lambda\in[0,1]$, 
	\[
	\Phi_k\big(\mathrm{e}^{\lambda x_k+(1-\lambda)y_k}, \cdots, \mathrm{e}^{\lambda x_L+(1-\lambda)y_L}\big)>0~(k=1,\cdots,L).
	\]
	Actually, we can show that
	\begin{equation}\label{Prop-ReductionF-Rienhardt-Proof-2}
	\Phi_k\big(\mathrm{e}^{\lambda x_k+(1-\lambda)y_k}, \cdots, \mathrm{e}^{\lambda x_L+(1-\lambda)y_L}\big)
	\geq \Phi_k\big(\mathrm{e}^{x_k}, \cdots, \mathrm{e}^{x_L}\big)^\lambda 
	\Phi_k\big(\mathrm{e}^{y_k}, \cdots, \mathrm{e}^{y_L}\big)^{1-\lambda}>0
	\end{equation}
	for $k=1,\cdots,L$. For brevity, we write 
	$X_k\equiv \Phi_k\big(\mathrm{e}^{x_k}, \cdots, \mathrm{e}^{x_L}\big)$, $Y_k\equiv \Phi_k\big(\mathrm{e}^{y_k}, \cdots, \mathrm{e}^{y_L}\big)$, and
	\[ 
	C_k\equiv \Phi_k\big(\mathrm{e}^{\lambda x_k+(1-\lambda)y_k}, \cdots, \mathrm{e}^{\lambda x_L+(1-\lambda)y_L}\big),
	\]
	so that \eqref{Prop-ReductionF-Rienhardt-Proof-2} can be written simply as $C_k\geq X_k^\lambda Y_k^{1-\lambda}>0$. When $k=L$, by the concavity of function $u\mapsto \log(1-\mathrm{e}^u)$, we have 
	\[
	\log\big(1-\mathrm{e}^{\lambda x_L+(1-\lambda)y_L}\big)
	\leq\lambda\log\big(1-\mathrm{e}^{x_L}\big)+(1-\lambda)\log\big(1-\mathrm{e}^{y_L}\big),
	\] 
	or equivalently, $C_L\geq X_L^\lambda Y_L^{1-\lambda}>0$. Next, assume that 
	$C_{k+1}\geq X_{k+1}^\lambda 
	Y_{k+1}^{1-\lambda}>0$.	Since $\mathrm{e}^{x_k}<X_{k+1}$ and $\mathrm{e}^{y_k}<Y_{k+1}$, we have
	$\mathrm{e}^{\lambda x_k+(1-\lambda)y_k}
	<X_{k+1}^\lambda Y_{k+1}^{1-\lambda}
	\leq C_{k+1}$. Thus
	\begin{align*}
	\log C_k
	=\log \left(1-\frac{\mathrm{e}^{\lambda x_k+(1-\lambda)y_k}}{C_{k+1}}\right)
	&\geq \log \left(1-\frac{\mathrm{e}^{\lambda x_k+(1-\lambda)y_k}}{X_{k+1}^\lambda 
		Y_{k+1}^{1-\lambda}}\right)\\
	&\geq \lambda\log\left(1-\frac{\mathrm{e}^{x_k}}{X_{k+1}}\right)
	+(1-\lambda)\log\left(1-\frac{\mathrm{e}^{y_k}}{Y_{k+1}}\right)\\
	&=\lambda \log X_k+(1-\lambda)\log Y_k.
	\end{align*}
	This proves \eqref{Prop-ReductionF-Rienhardt-Proof-2}. 
\end{proof}

\section{Euler-type integral representations of $F_{K}^{(L)}$}\label{Sect-3}

First of all, let us introduce the \emph{Dirichlet measure} $\mu_{\alpha,\beta}$ defined by \cite[p. 64, Definition 4.4-1]{Carlson-Book-1977}
\begin{equation}\label{DirichletM-1}
	\mathrm{d}\mu_{\alpha,\beta}(t):=\frac{\Gamma(\alpha+\beta)}{\Gamma(\alpha)\Gamma(\beta)}t^{\alpha-1}(1-t)^{\beta-1}\mathrm{d}t, \quad \min\{\Re(\alpha),\Re(\beta)\}>0. 
\end{equation}

Rosenhaus \cite[p. 24]{Rosenhaus-2019} has derived the following double integral representation for $F_K^{(L)}$ function. 
\begin{proposition}\label{Prop-RosenhausI}
	Let $L\geq4$ and $\mathbf{z}\in\mathbb{D}_L$. We have
	\begin{align}\label{IntegralR-FKL}
		&F_K^{(L)}\left[
		\begin{matrix}
			a_1, b_1,\cdots,b_{L-1}, a_2 \\ 
			c_1,\cdots,c_L
		\end{matrix}; \mathbf{z}
		\right]=\int_{0}^{1}\int_{0}^{1}
		(1-t_1z_1)^{-b_1}(1-t_Lz_L)^{-b_{L-1}}\notag\\
		&\hspace{0.5cm}\cdot F_K^{(L-2)}\left[
		\begin{matrix}
			b_1,\cdots,b_{L-1}\\
			c_2,\cdots,c_{L-1}
		\end{matrix};
		\frac{z_2}{1-t_1z_1},z_3,\cdots,z_{L-2},\frac{z_{L-1}}{1-t_Lz_L}
		\right]
		\mathrm{d}\mu_{a_1,c_1-a_1}(t_1)
		\mathrm{d}\mu_{a_2,c_L-a_2}(t_L),
	\end{align}
	where $\Re(c_1)>\Re(a_1)>0$, $\Re(c_L)>\Re(a_2)>0$. For $L=3$, the integral remains valid if the first and last variables are combined into a single variable.
\end{proposition}
\begin{proof}
	The integral itself can be easily established by expanding the integrand as a multivariate power series, interchanging the order of integration and summation, and then evaluating the beta integrals. Moreover, the conditions $\Re(c_1)>\Re(a_1)>0$, $\Re(c_L)>\Re(a_2)>0$ are easy to determine. However, Rosenhaus \cite[p. 24]{Rosenhaus-2019} did not specify the exact range of $\mathbf{z}$, which, in our opinion, requires some  explanation. To show that the integral formula \eqref{IntegralR-FKL} holds for $\mathbf{z}\in\mathbb{D}_L$, we prove that 
	\begin{equation}\label{Prop-RosenhausI-Proof-1}
	\mathbf{z}\in\mathbb{D}_L
	~\iff~
	(r_1<1)\wedge (r_L<1)\wedge
	\left(\left(\frac{r_2}{1-r_1},r_3,\cdots,r_{L-2},\frac{r_{L-1}}{1-r_L}\right)\in \mathbb{D}_{L-2}\right),
	\end{equation}
	where $r_j=|z_j|$ $(j=1,\cdots,L)$. We will use this result in Theorem \ref{Th-IntegralR} below. 
	
	Actually, if $\mathbf{z}\in\mathbb{D}_L$, we have immediately $r_L<1$ and $r_1<1$. Then, let 
	\[
	(s_1,\cdots,s_{L-2}):=\left(\frac{r_2}{1-r_1},r_3,\cdots,r_{L-2},\frac{r_{L-1}}{1-r_L}\right).
	\]
	From the definition of $\mathbb{D}_{L}$, we know that $(s_1,\cdots,s_{L-2})\in \mathbb{D}_{L-2}$ means that
	\[
	\Lambda_j(s_j,\cdots,s_{L-2})>0, \quad j=1,\cdots,L-2,
	\] 
	where $\Lambda_{L-2}(s_{L-2})=1-s_{L-2}$, and 
	\[
	\Lambda_j(s_j,\cdots,s_{L-2})=1-\frac{s_j}{\Lambda_{j+1}(s_{j+1},\cdots,s_{L-2})}, \quad j=1,\cdots,L-3. 
	\]
	But from $\mathbf{z}\in\mathbb{D}_L$, we have
	\begin{align*} 
	\Lambda_{L-2}(s_{L-2})
	&=1-\frac{r_{L-1}}{1-r_L}=\Phi_{L-1}(r_{L-1},r_L)>0,\\
	\Lambda_{j}(s_j,\cdots,s_{L-2})&=\Phi_{j+1}(r_{j+1},\cdots,r_L)>0, \quad j=2,\cdots, L-2,
	\end{align*}
	and
	\[
	\Lambda_1(s_1,\cdots,s_{L-2})=1-\frac{s_1}{\Lambda_2(s_2,\cdots,s_{L-2})}=1-\frac{r_2}{(1-r_1)\Phi_3(r_3,\cdots,r_L)}>0.
	\]
	Thus, $(s_1,\cdots,s_{L-2})\in \mathbb{D}_{L-2}$.
	
	Conversely, from $r_1<1$, $r_L<1$, and $(s_1,\cdots,s_{L-2})\in \mathbb{D}_{L-2}$, we have $\Lambda_j(s_j,\cdots,s_{L-2})>0$ $(j=1,\cdots,L-2)$. Then, based on the above derivation, we immediately obtain 
	\[
	\Phi_k(r_k,\cdots,r_L)=\Lambda_{k-1}(s_{k-1},\cdots,s_{L-2})>0, \quad k=3,\cdots, L-1.
	\]
	From $\Lambda_1(s_1,\cdots,s_{L-2})>0$, we have
	\[
	1-\frac{r_2}{(1-r_1)\Phi_3(r_3,\cdots,r_L)}>0
	\]
	and thus $0<r_1<\Phi_2(r_2,\cdots,r_L)$. Therefore,  $\Phi_1(r_1,\cdots,r_L)=1-\dfrac{r_1}{\Phi_2(r_2,\cdots,r_L)}>0$.
	Finally, $r_L<1$ guarantees that $\Phi_L(r_L)>0$. Hence, $\mathbf{z}\in\mathbb{D}_L$.
\end{proof}
\begin{remark}
	Usually, an Euler-type integral representation provides an analytic continuation of a multivariate hypergeometric function. However, integral \eqref{IntegralR-FKL} has functions of the same type nested in its integrand, so one cannot simply write down the domain of its analytic continuation.
\end{remark}

When $L=4$, the function $F_K^{(L-2)}$ on the right-hand side of \eqref{IntegralR-FKL} has only the first and last arguments, i.e., 
\begin{align*}
	&F_K^{(4)}\left[
	\begin{matrix}
		a_1, b_1,b_2,b_3, a_2 \\ 
		c_1,c_2,c_3,c_4
	\end{matrix}; z_1,z_2,z_3,z_4
	\right]=\int_{0}^{1}\int_{0}^{1}
	(1-t_1z_1)^{-b_1}(1-t_4z_4)^{-b_3}\notag\\
	&\hspace{0.5cm}\cdot F_K^{(2)}\left[
	\begin{matrix}
		b_1,b_2,b_3\\
		c_2,c_3
	\end{matrix};
	\frac{z_2}{1-t_1z_1},\frac{z_3}{1-t_4z_4}
	\right]
	\mathrm{d}\mu_{a_1,c_1-a_1}(t_1)
	\mathrm{d}\mu_{a_2,c_4-a_2}(t_4),
\end{align*}
where $\Re(c_1)>\Re(a_1)>0$, $\Re(c_4)>\Re(a_2)>0$ and $(z_1,z_2,z_3,z_4)\in\mathbb{D}_{4}$. 

When $L=3$, we have
\begin{align}\label{IntegralR-FKL-L3}
	&F_K^{(3)}\left[
	\begin{matrix}
		a_1, b_1,b_2, a_2 \\ 
		c_1,c_2,c_3
	\end{matrix}; z_1,z_2,z_3
	\right]=\int_{0}^{1}\int_{0}^{1}
	(1-t_1z_1)^{-b_1}(1-t_3z_3)^{-b_2}\notag\\
	&\hspace{1.5cm}\cdot{}_{2}F_{1}\left[
	\begin{matrix}
		b_1,b_2\\
		c_2
	\end{matrix};
	\frac{z_2}{(1-t_1z_1)(1-t_3z_3)}
	\right]
	\mathrm{d}\mu_{a_1,c_1-a_1}(t_1)
	\mathrm{d}\mu_{a_2,c_3-a_2}(t_3),
\end{align}
where $\Re(c_1)>\Re(a_1)>0$, $\Re(c_3)>\Re(a_2)>0$ and $(z_1,z_2,z_3)\in\mathbb{D}_3$. This formula provides a new integral representation for Saran's $F_K$ function defined by \eqref{Def-SaranFK}. Furthermore, since ${}_2F_1$ admits analytic continuation to $\mathbb{C}\setminus[1,+\infty)$, the integral on the right-hand side of \eqref{IntegralR-FKL-L3} provides an analytic continuation of $F_K^{(3)}$ to the set
\[
\widetilde{\mathbb{D}_3}:=\left\{(z_1,z_2,z_3)\in\mathbb{C}^3;~
z_1,z_3\notin[1,+\infty)~\text{and}~ \frac{z_2}{(1-t_1z_1)(1-t_3z_3)}\notin[1,+\infty)~
(\forall t_1,t_3\in[0,1])\right\}.
\]

If we further replace the Gauss hypergeometric ${}_2F_1$ in \eqref{IntegralR-FKL} by its Euler-type integral representation \cite[p. 77, Theorem 5.3]{Kristensson-Book-2010}
\begin{equation}\label{IntegralR-FKL-L1}
	{}_{2}F_{1}\left[\begin{matrix}
		a, b \\
		c
	\end{matrix};z\right]
	=\int_0^1 (1-zt)^{-a}\mathrm{d}\mu_{b,c-b}(t),\quad \Re(c)>\Re(b)>0, ~z\notin[1,+\infty),
\end{equation}
we reproduce Saran's triple integral \cite[p. 299]{Saran-1955a}: 
\begin{align}\label{IntegralR-FKL-L3-Saran}
	&F_K^{(3)}\left[
	\begin{matrix}
		a_1, b_1,b_2, a_2 \\ 
		c_1,c_2,c_3
	\end{matrix}; z_1,z_2,z_3
	\right]=\int_{(0,1)^3}
	(1-t_1z_1)^{b_2-b_1}\notag\\
	&\hspace{0.5cm}\cdot
	\left(1-t_1z_1-t_3z_3+t_1t_3z_1z_3-z_2t_2\right)^{-b_2}
	\mathrm{d}\mu_{a_1,c_1-a_1}(t_1)
	\mathrm{d}\mu_{b_1,c_2-b_1}(t_2)
	\mathrm{d}\mu_{a_2,c_3-a_2}(t_3),
\end{align}
where $\Re(c_1)>\Re(a_1)>0$, $\Re(c_2)>\Re(b_1)>0$, $\Re(c_3)>\Re(a_2)>0$, and $(z_1,z_2,z_3)\in\mathbb{D}_3$. 

In the remainder of this section, we shall point out that, by using the integral representation \eqref{IntegralR-FKL} repeatedly, we can eliminate the $F_K^{(L-2)}$ function from the integrand and obtain an $L$-fold integral whose integrand is an elementary function. 

Recall that Proposition \ref{Prop-DecomF} of Section \ref{Sect-2} introduces a system of functions $\Phi_k$ $(k=1,\cdots L)$. In the subsequent discussion, we need to define another system of functions \emph{in the reverse direction}. Let $\Psi_1(z_1)=1-z_1$ and  
\[
\Psi_k(z_1,\cdots,z_k)=1-\frac{z_k}{\Psi_{k-1}(z_1,\cdots,z_{k-1})}\quad (k=2,\cdots,L).
\]
As before, we write $\Psi_k\equiv \Psi_k(z_1,\cdots,z_k)$, provided that no confusion arises.
\begin{theorem}\label{Th-IntegralR}
	Let $L\geq3$ and $\mathbf{z}\in\mathbb{D}_L$. If $L=2M$, then we have
	\begin{equation}\label{Th-IntegralR-1}
		F_K^{(L)}\left[
		\begin{matrix}
			a_1, b_1,\cdots,b_{L-1}, a_2 \\ 
			c_1,\cdots,c_L
		\end{matrix}; \mathbf{z}
		\right]=\int_{(0,1)^{L}}
		\prod_{i=1}^{M}\widehat{\Psi}_i^{-b_i}
		\prod_{j=M+1}^{L}\widehat{\Phi}_j^{-b_{j-1}}
		\prod_{\ell=1}^{L}\mathrm{d}\mu_{\alpha_\ell,\gamma_\ell-\alpha_\ell}(t_\ell),
	\end{equation}
	where  
	\[
	(\alpha_\ell,\gamma_\ell)=\begin{cases}
		(a_1,c_1), & \ell=1, \\
		(b_{\ell-1},c_\ell), & 2\leq \ell\leq M, \\
		(b_\ell,c_\ell), & M+1\leq \ell\leq L-1, \\
		(a_2,c_L), & \ell=L,
	\end{cases}
	\] 
	and $\Re(\gamma_\ell)>\Re(\alpha_\ell)>0$ $(\ell=1,\cdots,L)$. If
	$L=2M+1$, then we have  
	\begin{align}\label{Th-IntegralR-2}
		F_K^{(L)}\left[
		\begin{matrix}
			a_1, b_1,\cdots,b_{L-1}, a_2 \\ 
			c_1,\cdots,c_L
		\end{matrix}; \mathbf{z}
		\right]&=\int_{(0,1)^L}
		\prod_{i=1}^{M}\widehat{\Psi}_i^{-b_i}
		\prod_{j=M+2}^{L}\widehat{\Phi}_j^{-b_{j-1}}\notag\\
		&\hspace{0.5cm}\cdot 
		\left(1-\frac{z_{M+1}t_{M+1}}{\widehat{\Psi}_M\widehat{\Phi}_{M+2}}\right)^{-b_{M+1}}
		\prod_{\ell=1}^{L}\mathrm{d}\mu_{\alpha_\ell,\gamma_\ell-\alpha_\ell}(t_\ell),
	\end{align}
	where  
	\[
	(\alpha_\ell,\gamma_\ell)=\begin{cases}
		(a_1,c_1), & \ell=1, \\
		(b_{\ell-1},c_\ell), & 2\leq \ell\leq M+1, \\
		(b_\ell,c_\ell), & M+2\leq \ell\leq L-1, \\
		(a_2,c_L), & \ell=L,
	\end{cases}
	\] 
	and $\Re(\gamma_\ell)>\Re(\alpha_\ell)>0$ $(\ell=1,\cdots,L)$.
	Here, we write $\widehat{\Psi}_j\equiv\Psi(t_1z_1,\cdots,t_jz_j)$ and $\widehat{\Phi}_j\equiv\Phi(t_jz_j,\cdots,t_Lz_L)$, for short. 
\end{theorem}
\begin{proof}
	Let us apply the same integral representation \eqref{IntegralR-FKL} to the inner function 
	\[
	F_K^{(L-2)}\left[
	\begin{matrix}
		b_1,\cdots,b_{L-1}\\
		c_2,\cdots,c_{L-1}
	\end{matrix};
	w_1,\cdots,w_{L-2}
	\right],
	\]
	where $w_1\equiv\dfrac{z_2}{1-t_1z_1}$, $w_2\equiv z_3, \cdots, w_{L-3}\equiv z_{L-2}$, and $w_{L-2}\equiv \dfrac{z_{L-1}}{1-t_Lz_L}$. Then we obtain 
	\begin{align}\label{Th-IntegralR-Proof-1}
		&F_K^{(L-2)}\left[
		\begin{matrix}
			b_1,\cdots,b_{L-1}\\
			c_2,\cdots,c_{L-1}
		\end{matrix};
		w_1,\cdots,w_{L-2}
		\right]
		=\int_{0}^{1}\int_{0}^{1}
		(1-t_2w_1)^{-b_2}(1-t_{L-1}w_{L-2})^{-b_{L-2}}\notag\\
		&\hspace{0.5cm}\cdot F_K^{(L-4)}\left[
		\begin{matrix}
			b_2,\cdots,b_{L-2}\\
			c_3,\cdots,c_{L-2}
		\end{matrix};
		\frac{w_2}{1-t_2w_1},w_3,\cdots,w_{L-4},\frac{w_{L-3}}{1-t_{L-1}w_{L-2}}
		\right]\notag\\
		&\hspace{0.5cm}\cdot 
		\mathrm{d}\mu_{b_1,c_2-b_1}(t_2)
		\mathrm{d}\mu_{b_{L-1},c_{L-1}-b_{L-1}}(t_{L-1}),
	\end{align}
	where $\Re(c_2)>\Re(b_1)>0$,  $\Re(c_{L-1})>\Re(b_{L-1})>0$, and 
	\begin{equation}\label{Th-IntegralR-Proof-2}
		\left(\frac{r_2}{1-r_1},r_3,\cdots,r_{L-2},\frac{r_{L-1}}{1-r_L}\right)\in \mathbb{D}_{L-2}.
	\end{equation} 
	Here, as usual, $r_j=|z_j|$ $(j=1,\cdots,L)$. In view of \eqref{Prop-RosenhausI-Proof-1}, it suffices to require that $\mathbf{z}\in\mathbb{D}_L$, without imposing the additional assumption \eqref{Th-IntegralR-Proof-2}. Note that
	\begin{align}\label{Th-IntegralR-Proof-3}
		&(1-t_1z_1)^{-b_1}(1-t_2w_1)^{-b_2}
		(1-t_{L-1}w_{L-2})^{-b_{L-2}}(1-t_Lz_L)^{-b_{L-1}}\notag\\
		&\hspace{0.5cm}=(1-t_1z_1)^{-b_1}\left(1-\frac{t_2z_2}{1-t_1z_1}\right)^{-b_2}\left(1-\frac{t_{L-1}z_{L-1}}{1-t_Lz_L}\right)^{-b_{L-2}}(1-t_L z_L)^{-b_{L-1}}\notag\\
		&\hspace{0.5cm}=\Psi_1^{-b_1}(t_1z_1)
		\Psi_2^{-b_2}(t_1z_1,t_2z_2)
		\Phi_{L-1}^{-b_{L-2}}(t_{L-1}z_{L-1},t_Lz_L)
		\Phi_L^{-b_{L-1}}(t_Lz_L)\notag\\
		&\hspace{0.5cm}=\widehat{\Psi}_1^{-b_1}
		\widehat{\Psi}_2^{-b_2}
		\widehat{\Phi}_{L-1}^{-b_{L-2}}\widehat{\Phi}_L^{-b_{L-1}}.
	\end{align}
	Thus, substituting \eqref{Th-IntegralR-Proof-1} into \eqref{IntegralR-FKL} and taking into account \eqref{Th-IntegralR-Proof-3}, we have
	\begin{align}\label{Th-IntegralR-Proof-4}
		&F_K^{(L)}\left[
		\begin{matrix}
			a_1, b_1,\cdots,b_{L-1}, a_2 \\ 
			c_1,\cdots,c_L
		\end{matrix}; \mathbf{z}
		\right]=\int_{(0,1)^{4}}
		\widehat{\Psi}_1^{-b_1}
		\widehat{\Psi}_2^{-b_2}
		\widehat{\Phi}_{L-1}^{-b_{L-2}}
		\widehat{\Phi}_L^{-b_{L-1}}\notag\\
		&\hspace{1cm}\cdot F_K^{(L-4)}\left[
		\begin{matrix}
			b_2,\cdots,b_{L-2}\\
			c_3,\cdots,c_{L-2}
		\end{matrix};
		\frac{z_3}{\widehat{\Psi}_2},z_4,\cdots,z_{L-3},\frac{z_{L-2}}{\widehat{\Phi}_{L-1}}
		\right]\notag\\
		&\hspace{1cm}\cdot\mathrm{d}\mu_{a_1,c_1-a_1}(t_1)
		\mathrm{d}\mu_{b_1,c_2-b_1}(t_2)
		\mathrm{d}\mu_{b_{L-1},c_{L-1}-b_{L-1}}(t_{L-1})
		\mathrm{d}\mu_{a_2,c_L-a_2}(t_L),
	\end{align}
	where $\Re(c_1)>\Re(a_1)>0$, $\Re(c_2)>\Re(b_1)>0$, $\Re(c_{L-1})>\Re(b_{L-1})>0$, $\Re(c_L)>\Re(a_2)>0$, and $\mathbf{z}\in\mathbb{D}_{L}$.

	Similarly, we consider
	\[
	F_K^{(L-4)}\left[
	\begin{matrix}
		b_2,\cdots,b_{L-2}\\
		c_3,\cdots,c_{L-2}
	\end{matrix};
	u_1,\cdots, u_{L-4}
	\right]
	\]
	where $u_1\equiv\dfrac{z_3}{\widehat{\Psi}_2}, u_2\equiv z_4, \cdots, u_{L-4}\equiv \dfrac{z_{L-2}}{\widehat{\Phi}_{L-1}}$. Its integral representation is
	\begin{align*}
		&F_K^{(L-4)}\left[
		\begin{matrix}
			b_2,\cdots,b_{L-2}\\
			c_3,\cdots,c_{L-2}
		\end{matrix};
		u_1,\cdots,u_{L-4}
		\right]
		=\int_{0}^{1}\int_{0}^{1}
		(1-t_3u_1)^{-b_3}(1-t_{L-2}u_{L-4})^{-b_{L-3}}\notag\\
		&\hspace{0.5cm}\cdot F_K^{(L-6)}\left[
		\begin{matrix}
			b_3,\cdots,b_{L-3}\\
			c_4,\cdots,c_{L-3}
		\end{matrix};
		\frac{u_2}{\widehat{\Psi}_3},u_3,\cdots,u_{L-6},\frac{u_{L-5}}{\widehat{\Phi}_{L-2}}
		\right]\mathrm{d}\mu_{b_2,c_3-b_2}(t_3)
		\mathrm{d}\mu_{b_{L-2},c_{L-2}-b_{L-2}}(t_{L-2}),
	\end{align*}
	where $\Re(c_3)>\Re(b_2)>0$, $\Re(c_{L-2})>\Re(b_{L-2})>0$, and $\mathbf{z}\in\mathbb{D}_L$. Substituting this expressions into \eqref{Th-IntegralR-Proof-4} and taking into account that $(1-t_3u_1)^{-b_3}(1-t_{L-2}u_{L-4})^{-b_{L-3}}
	=\widehat{\Psi}_3^{-b_3}
	\widehat{\Phi}_{L-2}^{-b_{L-3}}$, we have 
	\begin{align*}
		&F_K^{(L)}\left[
		\begin{matrix}
			a_1, b_1,\cdots,b_{L-1}, a_2 \\ 
			c_1,\cdots,c_L
		\end{matrix}; \mathbf{z}
		\right]=\int_{(0,1)^{6}}
		\widehat{\Psi}_1^{-b_1}
		\widehat{\Psi}_2^{-b_2}
		\widehat{\Psi}_3^{-b_3}
		\widehat{\Phi}_{L-2}^{-b_{L-3}}
		\widehat{\Phi}_{L-1}^{-b_{L-2}}
		\widehat{\Phi}_L^{-b_{L-1}}\notag\\
		&\hspace{1cm}\cdot F_K^{(L-6)}\left[
		\begin{matrix}
			b_3,\cdots,b_{L-3}\\
			c_4,\cdots,c_{L-3}
		\end{matrix};
		\frac{z_4}{\widehat{\Psi}_3},z_5,\cdots,z_{L-4},\frac{z_{L-3}}{\widehat{\Phi}_{L-2}}
		\right]\mathrm{d}\mu_{a_1,c_1-a_1}(t_1)
		\mathrm{d}\mu_{b_1,c_2-b_1}(t_2)\notag\\
		&\hspace{1cm}\cdot \mathrm{d}\mu_{b_2,c_3-b_2}(t_3)
		\mathrm{d}\mu_{b_{L-2},c_{L-2}-b_{L-2}}(t_{L-2})
		\mathrm{d}\mu_{b_{L-1},c_{L-1}-b_{L-1}}(t_{L-1})
		\mathrm{d}\mu_{a_2,c_L-a_2}(t_L).
	\end{align*}
	where $\Re(c_1)>\Re(a_1)>0$, 
	\begin{align*}
	\Re(c_{\ell})>\Re(b_{\ell-1})>0 \quad &(\ell=2,3),\\
	\Re(c_{\ell})>\Re(b_{\ell})>0 \quad &(\ell=L-2,L-1), 
	\end{align*}
	$\Re(c_L)>\Re(a_2)>0$, and $\mathbf{z}\in\mathbb{D}_L$. The pattern is now clear. 
	
	Assume that after $m$ substitutions $(0\leq m<\lfloor L/2\rfloor)$, we have
	\begin{align*}
		&F_K^{(L)}\left[
		\begin{matrix}
			a_1, b_1,\cdots,b_{L-1}, a_2 \\ 
			c_1,\cdots,c_L
		\end{matrix}; \mathbf{z}
		\right]=\int_{(0,1)^{2m}}
		\prod_{i=1}^{m}\widehat{\Psi}_i^{-b_i}
		\prod_{j=0}^{m-1}\widehat{\Phi}_{L-j}^{-b_{L-1-j}}\notag\\
		&\hspace{0.5cm}\cdot F_K^{(L-2m)}\left[
		\begin{matrix}
			b_m,\cdots,b_{L-m}\\
			c_{m+1},\cdots,c_{L-m}
		\end{matrix};
		\frac{z_{m+1}}{\widehat{\Psi}_m},z_{m+2},\cdots,z_{L-m-1},\frac{z_{L-m}}{\widehat{\Phi}_{L-m+1}}
		\right]\notag\\
		&\hspace{0.5cm}\cdot\prod_{i=1}^{m}\mathrm{d}\mu_{\alpha_i,\gamma_i-\alpha_i}(t_i)
		\prod_{j=0}^{m-1}\mathrm{d}\mu_{\alpha_{L-j},\gamma_{L-j}-\alpha_{L-j}}(t_{L-j}).
	\end{align*}
	Applying the recursive integral representation once more to the remaining $F_K^{(L-2m)}$ function, we obtain
	\begin{align*}
		&F_K^{(L)}\left[
		\begin{matrix}
			a_1, b_1,\cdots,b_{L-1}, a_2 \\ 
			c_1,\cdots,c_L
		\end{matrix}; \mathbf{z}
		\right]=\int_{(0,1)^{2m+2}}
		\prod_{i=1}^{m+1}\widehat{\Psi}_i^{-b_i}
		\prod_{j=0}^{m}\widehat{\Phi}_{L-j}^{-b_{L-1-j}}\notag\\
		&\hspace{0.5cm}\cdot F_K^{(L-2m-2)}\left[
		\begin{matrix}
			b_{m+1},\cdots,b_{L-m-1}\\
			c_{m+2},\cdots,c_{L-m-1}
		\end{matrix};
		\frac{z_{m+2}}{\widehat{\Psi}_{m+1}},z_{m+3},\cdots,z_{L-m-2},\frac{z_{L-m-1}}{\widehat{\Phi}_{L-m}}
		\right]\notag\\
		&\hspace{0.5cm}\cdot\prod_{i=1}^{m+1}\mathrm{d}\mu_{\alpha_i,\gamma_i-\alpha_i}(t_i)
		\prod_{j=0}^{m}\mathrm{d}\mu_{\alpha_{L-j},\gamma_{L-j}-\alpha_{L-j}}(t_{L-j}). 
	\end{align*}
	This is exactly the induction hypothesis with $m$ replaced by $m+1$.
	
	After $\lfloor L/2\rfloor$ steps, the $F_K^{(L)}$ function disappears entirely when $L$ is even, or reduces to a Gauss hypergeometric function ${}_{2}F_{1}$ when $L$ is odd. In the first case, we obtain 
	\begin{align*}
		F_K^{(L)}\left[
		\begin{matrix}
			a_1, b_1,\cdots,b_{L-1}, a_2 \\ 
			c_1,\cdots,c_L
		\end{matrix}; \mathbf{z}
		\right]
		&=\int_{(0,1)^L}
		\prod_{i=1}^{M}\widehat{\Psi}_i^{-b_i}
		\prod_{j=0}^{M-1}\widehat{\Phi}_{L-j}^{-b_{L-1-j}}\notag\\
		&\hspace{1.5cm}\cdot\prod_{i=1}^{M}\mathrm{d}\mu_{\alpha_i,\gamma_i-\alpha_i}(t_i)
		\prod_{j=0}^{M-1}\mathrm{d}\mu_{\alpha_{L-j},\gamma_{L-j}-\alpha_{L-j}}(t_{L-j})\notag\\
		&=\int_{(0,1)^L}
		\prod_{i=1}^{M}\widehat{\Psi}_i^{-b_i}
		\prod_{j=M+1}^{L}\widehat{\Phi}_j^{-b_{j-1}}\prod_{\ell=1}^{L}\mathrm{d}\mu_{\alpha_\ell,\gamma_\ell-\alpha_\ell}(t_\ell),
	\end{align*}
	where $M=L/2$. In the second case, we have 
	\begin{align}\label{Th-IntegralR-Proof-5} 
		&F_K^{(L)}\left[
		\begin{matrix}
			a_1, b_1,\cdots,b_{L-1}, a_2 \\ 
			c_1,\cdots,c_L
		\end{matrix}; \mathbf{z}
		\right]=\int_{(0,1)^{L-1}}
		\prod_{i=1}^{M}\widehat{\Psi}_i^{-b_i}
		\prod_{j=0}^{M-1}\widehat{\Phi}_{L-j}^{-b_{L-1-j}}\notag\\
		&\hspace{0.5cm}\cdot {}_{2}F_{1}\left[
		\begin{matrix}
			b_M,b_{M+1}\\
			c_{M+1}
		\end{matrix};
		\frac{z_{M+1}}{\widehat{\Psi}_M\widehat{\Phi}_{M+2}}
		\right]\prod_{i=1}^{M}\mathrm{d}\mu_{\alpha_i,\gamma_i-\alpha_i}(t_i)
		\prod_{j=0}^{M-1}\mathrm{d}\mu_{\alpha_{L-j},\gamma_{L-j}-\alpha_{L-j}}(t_{L-j}),
	\end{align}
	where $M=(L-1)/2$. By \eqref{IntegralR-FKL-L1}, we can further express \eqref{Th-IntegralR-Proof-5} as follows: 
	\begin{align*}
		F_K^{(L)}\left[
		\begin{matrix}
			a_1, b_1,\cdots,b_{L-1}, a_2 \\ 
			c_1,\cdots,c_L
		\end{matrix}; \mathbf{z}
		\right]&=\int_{(0,1)^L}
		\prod_{i=1}^{M}\widehat{\Psi}_i^{-b_i}\prod_{j=M+2}^{L}\widehat{\Phi}_j^{-b_{j-1}}\notag\\
		&\hspace{0.5cm}\cdot 
		\left(1-\frac{z_{M+1}t_{M+1}}{\widehat{\Psi}_M\widehat{\Phi}_{M+2}}\right)^{-b_{M+1}}
		\prod_{\ell=1}^{L}\mathrm{d}\mu_{\alpha_\ell,\gamma_\ell-\alpha_\ell}(t_\ell).
	\end{align*}
The proof is therefore complete.
\end{proof}
\begin{remark}
When $L=3$, the integral \eqref{Th-IntegralR-2} becomes
\begin{align}\label{Th-IntegralR-Remark-1}
	F_K^{(3)}\left[
	\begin{matrix}
		a_1, b_1, b_2, a_2 \\ 
		c_1, c_2, c_3
	\end{matrix}; z_1,z_2,z_3
	\right]=\int_{(0,1)^3}
	\widehat{\Psi}_1^{b_2-b_1}
	\left(\widehat{\Psi}_1\widehat{\Phi}_{3}-z_2t_2\right)^{-b_2}
	\prod_{\ell=1}^3\mathrm{d}\mu_{\alpha_\ell,\gamma_\ell-\alpha_\ell}(t_\ell),
\end{align}  
where
\begin{align}
	\widehat{\Psi}_1&=1-t_1z_1,\label{Th-IntegralR-Remark-2}\\
	\widehat{\Phi}_3&=1-t_3z_3,\label{Th-IntegralR-Remark-3}
\end{align}
and
\begin{equation}\label{Th-IntegralR-Remark-4}
(\alpha_\ell,\gamma_\ell)=\begin{cases}
	(a_1,c_1), & \ell=1, \\
	(b_1,c_2), & \ell=2, \\
	(a_2,c_3), & \ell=3.
\end{cases}
\end{equation}
Substituting \eqref{Th-IntegralR-Remark-2}--\eqref{Th-IntegralR-Remark-4} into \eqref{Th-IntegralR-Remark-1}, we obtain the classical result \eqref{IntegralR-FKL-L3-Saran}. When $L=4$, \eqref{Th-IntegralR-1} reduces to
\begin{equation}\label{Th-IntegralR-Remark-5}
	F_K^{\textcolor{orange}{(4)}}\left[
	\begin{matrix}
		a_1, b_1,b_2,b_3, a_2 \\ 
		c_1,c_2,c_3,c_4
	\end{matrix}; z_1,z_2,z_3,z_4
	\right]=\int_{(0,1)^4}
	\widehat{\Psi}_1^{-b_1}\widehat{\Psi}_2^{-b_2}
	\widehat{\Phi}_3^{-b_2}
	\widehat{\Phi}_4^{-b_3}
	\prod_{\ell=1}^{4}\mathrm{d}\mu_{\alpha_\ell,\gamma_\ell-\alpha_\ell}(t_\ell),
\end{equation}
where  
\begin{equation}\label{Th-IntegralR-Remark-6}
\widehat{\Psi}_1^{-b_1}\widehat{\Psi}_2^{-b_2}
\widehat{\Phi}_3^{-b_2}
\widehat{\Phi}_4^{-b_3}=(1-t_1z_1)^{b_2-b_1}
(1-t_1z_1-t_2z_2)^{-b_2}
(1-t_4z_4-t_3z_3)^{-b_2}
(1-t_4z_4)^{b_2-b_3}
\end{equation}
and
\begin{equation}\label{Th-IntegralR-Remark-7}
(\alpha_\ell,\gamma_\ell)=\begin{cases}
	(a_1,c_1), & \ell=1, \\
	(b_1,c_2), & \ell=2, \\
	(b_3,c_3), & \ell=3, \\
	(a_2,c_4), & \ell=4.
\end{cases}
\end{equation}
Substituting \eqref{Th-IntegralR-Remark-6} and \eqref{Th-IntegralR-Remark-7} into \eqref{Th-IntegralR-Remark-5} gives
\begin{align*}
	&F_K^{(4)}\left[
	\begin{matrix}
		a_1,b_1,b_2,b_3,a_2 \\ 
		c_1,c_2,c_3,c_4
	\end{matrix};z_1,z_2,z_3,z_4
	\right]
	=\int_{(0,1)^4}
	(1-t_1z_1)^{b_2-b_1}(1-t_4z_4)^{b_2-b_3}\notag\\
	&\hspace{2cm}\cdot\big(1-t_1z_1-t_2z_2-t_3z_3-t_4z_4
	+t_1t_3z_1z_3+t_1t_4z_1z_4+t_2t_3z_2z_3+t_2t_4z_2z_4\big)^{-b_2}\notag\\
	&\hspace{2cm}\cdot\mathrm{d}\mu_{a_1,c_1-a_1}(t_1)
	\mathrm{d}\mu_{b_1,c_2-b_1}(t_2)
	\mathrm{d}\mu_{b_3,c_3-b_3}(t_3)
	\mathrm{d}\mu_{a_2,c_4-a_2}(t_4),
\end{align*}
where $\Re(c_1)>\Re(a_1)>0$, $\Re(c_2)>\Re(b_1)>0$, $\Re(c_3)>\Re(b_3)>0$, $\Re(c_4)>\Re(a_2)>0$, and $(z_1,z_2,z_3,z_4)\in\mathbb{D}_4$.
\end{remark}

\section{Decomposition formulas}\label{Sect-4}

Let 
\begin{equation}\label{Def-SymbolicO-1}
\nabla_{z_i,z_j}(h)\equiv\frac{\Gamma(h)\Gamma(\delta_i+\delta_j+h)}{\Gamma(\delta _i+h)\Gamma(\delta_j+h)}
\end{equation}
and
\begin{equation}\label{Def-SymbolicO-2}
\Delta_{z_i,z_j}(h)\equiv\frac{\Gamma(\delta_i+h)\Gamma(\delta_j+h)}{\Gamma(h)\Gamma(\delta_i+\delta_j+h)},
\end{equation}
where $\delta_i\equiv z_i\dfrac{\partial}{\partial z_i}$, $\delta_j=z_j\dfrac{\partial}{\partial z_j}$ and $i\neq j$. These symbolic operators were introduced and studied by Burchnall and Chaundy (see \cite{Burchnall-Chaundy-1940} and \cite{Burchnall-Chaundy-1941}). It is easy to see that \cite[p. 250]{Burchnall-Chaundy-1940}
\begin{align*}
\nabla_{z_i,z_j}(h)\big\{(h)_m(h)_n z_i^m z_j^n\big\}
&=(h)_{m+n} z_i^m z_j^n,\\
\Delta_{z_i,z_j}(h)\big\{(h)_{m+n} z_i^m z_j^n\big\}
&=(h)_m(h)_n z_i^m z_j^n.
\end{align*}
On the other hand, due to the Gauss summation theorem \cite[p. 78, Corollary 5.1]{Kristensson-Book-2010}
\begin{equation}\label{Th-GaussianST}
	{}_{2}F_{1}\left[\begin{matrix}
		a,b\\
		c
	\end{matrix};1\right]=
	\frac{\Gamma(c)\Gamma(c-a-b)}{\Gamma(c-a)\Gamma(c-b)}, \quad \Re(c-a-b)>0, \quad c\in\mathbb{C}\setminus\mathbb{Z}_{\leq0}, 
\end{equation}
the operators \eqref{Def-SymbolicO-1} and \eqref{Def-SymbolicO-2} can also be expressed as
\begin{equation}\label{Def-SymbolicO-1a}
	\nabla_{z_i,z_j}(h)=\sum_{n=0}^{\infty}\frac{(-\delta_i)_n(-\delta_j)_n}{(h)_n~n!}
\end{equation}
and
\[
	\Delta_{z_i,z_j}(h)=\sum_{n=0}^{\infty}\frac{(-\delta_i)_n(-\delta_j)_n}{(1-h-\delta_i-\delta_j)_n~n!},
\]
respectively. If $a=-n$ in \eqref{Th-GaussianST}, it reduces to the Chu-Vandermonde formula \cite[p. 30, Eq. (8)]{Srivastava-Manocha-Book-1984}: 
\begin{equation}\label{Th-ChuVandermonde}
	{}_{2}F_{1}\left[\begin{matrix}
		-n,b\\
		c
	\end{matrix};1\right]=
	\frac{(c-b)_n}{(c)_n}, \quad c\in\mathbb{C}\setminus\mathbb{Z}_{\leq0}. 
\end{equation}

For every analytic function $f(z_i)$, the following operator identities \cite[p. 5, Eqs. (3.1) and (3.2)]{Hasanov-et al-2010}
\begin{equation}\label{PooleF-1}
	(-\delta_i)_n\big\{f(z_i)\big\}=(-1)^n z_i^n \frac{\mathrm{d}^n}{\mathrm{d}z_i^n}\big\{f(z_i)\big\}
\end{equation}
and 
\[
	(a+\delta_i)_n\big\{f(z_i)\big\}
	=z_i^{1-a}\frac{\mathrm{d}^n}{\mathrm{d}z_i^n}
	\big\{z^{a+n-1}f(z_i)\big\}
\]
are valid.

\begin{proposition}\label{Prop-DecomF}
	We have
	\begin{equation}\label{Prop-DecomF-1}
		F_K^{(L)}\left[
		\begin{matrix}
			a_1, b_1,\cdots,b_{L-1}, a_2 \\ 
			c_1,\cdots,c_L
		\end{matrix}; \mathbf{z}
		\right]=\nabla_{z_{L-1},z_{L}}(b_{L-1})\cdots \nabla_{z_1,z_2}(b_1)~\mathcal{F}(\mathbf{z})\textcolor{orange}{,}
	\end{equation}
	where
	\[
	\mathcal{F}(\mathbf{z}):=
	{}_2 F_1\left[\begin{matrix}
		a_2, b_{L-1}\\
		c_L
	\end{matrix};z_L\right]
	{}_2 F_1\left[\begin{matrix}
		b_1, a_1\\
		c_1
	\end{matrix};z_1\right]
	\prod_{k=2}^{L-1}{}_2 F_1\left[\begin{matrix}
		b_k, b_{k-1}\\
		c_k
	\end{matrix};z_k\right].
	\]
\end{proposition}
\begin{proof}
	Since
	\begin{align*}
	&\Delta_{z_{L-1},z_{L}}(b_{L-1})F_K^{(L)}\left[
	\begin{matrix}
		a_1, b_1,\cdots,b_{L-1}, a_2 \\ 
		c_1,\cdots,c_L
	\end{matrix}; \mathbf{z}
	\right]\\
	&\hspace{0.5cm}=\sum_{\mathbf{n}=\mathbf{0}}^{\infty} \frac{(a_1)_{n_1}(b_1)_{n_1+n_2}\cdots(b_{L-2})_{n_{L-2}+n_{L-1}}(a_2)_{n_L}}{(c_1)_{n_1}\cdots(c_L)_{n_L}} \frac{z_1^{n_1}}{n_1!} \cdots \frac{z_L^{n_{L-2}}}{n_{L-2}!}\\
	&\hspace{1.5cm}\cdot \Delta_{z_{L-1},z_{L}}(b_{L-1})\big[(b_{L-1})_{n_{L-1}+n_L} z_{L-1}^{n_{L-1}}z_L^{n_L}\big]\frac{1}{n_{L-1}!}\frac{1}{n_L!}\\
	&\hspace{0.5cm}=\sum_{\mathbf{n}=\mathbf{0}}^{\infty} \frac{(a_1)_{n_1}(b_1)_{n_1+n_2}\cdots(b_{L-2})_{n_{L-2}+n_{L-1}}(b_{L-1})_{n_{L-1}}}{(c_1)_{n_1}\cdots(c_{L-1})_{n_{L-1}}} 
	\frac{(a_2)_{n_L}(b_{L-1})_{n_L}}{(c_L)_{n_L}}\frac{\mathbf{z}^{\mathbf{n}}}{\mathbf{n}!}\\
	&\hspace{0.5cm}={}_2 F_1\left[\begin{matrix}
		a_2, b_{L-1}\\
		c_L
	\end{matrix};z_L\right] F_K^{(L-1)}\left[
	\begin{matrix}
		a_1, b_1,\cdots,b_{L-2}, b_{L-1} \\ 
		c_1,\cdots,c_{L-1}
	\end{matrix}; z_1,\cdots,z_{L-1}
	\right],
	\end{align*}
	we have
	\begin{align}\label{Prop-DecomF-Proof-1}
	&\prod_{k=1}^{L-1} \Delta_{z_{k},z_{k+1}}(b_k) 
	F_K^{(L)}\left[
	\begin{matrix}
		a_1, b_1,\cdots,b_{L-1}, a_2 \\ 
		c_1,\cdots,c_L
	\end{matrix}; \mathbf{z}
	\right]={}_2 F_1\left[\begin{matrix}
		a_2, b_{L-1}\\
		c_L
	\end{matrix};z_L\right]\notag\\
	&\hspace{0.5cm}\cdot \prod_{k=1}^{L-2} \Delta_{z_{k},z_{k+1}}(b_k) F_K^{(L-1)}\left[
	\begin{matrix}
		a_1, b_1,\cdots,b_{L-2}, b_{L-1} \\ 
		c_1,\cdots,c_{L-1}
	\end{matrix}; z_1,\cdots,z_{L-1}
	\right].
	\end{align}
	In similar manner, we have
	\begin{align}\label{Prop-DecomF-Proof-2}
		&\prod_{k=1}^{L-2} \Delta_{z_{k},z_{k+1}}(b_k) F_K^{(L-1)}\left[
		\begin{matrix}
			a_1, b_1,\cdots,b_{L-2}, b_{L-1} \\ 
			c_1,\cdots,c_{L-1}
		\end{matrix}; z_1,\cdots,z_{L-1}
		\right]\notag\\
		&\hspace{0.5cm}={}_2 F_1\left[\begin{matrix}
			b_{L-1}, b_{L-2}\\
			c_{L-1}
		\end{matrix};z_{L-1}\right] 
		\prod_{k=1}^{L-3} \Delta_{z_{k},z_{k+1}}(b_k) F_K^{(L-2)}\left[
		\begin{matrix}
			a_1, b_1,\cdots,b_{L-3}, b_{L-2} \\ 
			c_1,\cdots,c_{L-2}
		\end{matrix}; z_1,\cdots,z_{L-2}
		\right]\notag\\
		&\hspace{0.5cm}=\cdots={}_2 F_1\left[\begin{matrix}
			b_{L-1}, b_{L-2}\\
			c_{L-1}
		\end{matrix};z_{L-1}\right] 
		\cdots {}_2 F_1\left[\begin{matrix}
			b_3, b_2\\
			c_3
		\end{matrix};z_3\right]
		\Delta_{z_1,z_2}(b_1) F_K^{(2)}\left[
		\begin{matrix}
			a_1, b_1, b_2 \\ 
			c_1, c_2
		\end{matrix}; z_1,z_2\right],
	\end{align}
	and
	\begin{equation}\label{Prop-DecomF-Proof-3}
	\Delta_{z_1,z_2}(b_1) F_K^{(2)}\left[
	\begin{matrix}
		a_1, b_1, b_2 \\ 
		c_1, c_2
	\end{matrix}; z_1,z_2\right]
	={}_2 F_1\left[\begin{matrix}
		b_2, b_1\\
		c_2
	\end{matrix};z_2\right]
	{}_2 F_1\left[\begin{matrix}
		b_1, a_1\\
		c_1
	\end{matrix};z_1\right].
	\end{equation}
	Combining \eqref{Prop-DecomF-Proof-1}, \eqref{Prop-DecomF-Proof-2} and \eqref{Prop-DecomF-Proof-3} gives
	\begin{align*}
	&\prod_{k=1}^{L-1} \Delta_{z_{k},z_{k+1}}(b_k) 
	F_K^{(L)}\left[
	\begin{matrix}
		a_1, b_1,\cdots,b_{L-1}, a_2 \\ 
		c_1,\cdots,c_L
	\end{matrix}; \mathbf{z}
	\right]\\
	&\hspace{0.5cm}={}_2 F_1\left[\begin{matrix}
		a_2, b_{L-1}\\
		c_L
	\end{matrix};z_L\right]
	{}_2 F_1\left[\begin{matrix}
		b_1, a_1\\
		c_1
	\end{matrix};z_1\right]
	\prod_{k=2}^{L-1}{}_2 F_1\left[\begin{matrix}
		b_k, b_{k-1}\\
		c_k
	\end{matrix};z_k\right].
	\end{align*}
	taking into account the relationship between 
	$\nabla_{z_i,z_j}$ and $\Delta_{z_i,z_j}$, we immediately obtain \eqref{Prop-DecomF-1}. This completes the proof. 
\end{proof} 
\begin{remark}
	When $L=2$, \eqref{Prop-DecomF-1} reduces to 
	\begin{equation}\label{Prop-DecomF-Remark-1}
		F_K^{(2)}\left[\begin{matrix}
			a_1, b_1, a_2 \\ 
			c_1, c_2
		\end{matrix}; z_1,z_2\right]=\nabla_{z_1,z_2}(b_1)~{}_2 F_1\left[\begin{matrix}
			a_2, b_1\\
			c_2
		\end{matrix};z_2\right]
		{}_2 F_1\left[\begin{matrix}
			b_1, a_1\\
			c_1
		\end{matrix};z_1\right],
	\end{equation}
	which is the result of Burchnall and Chaundy \cite[p. 250, Eq. (6)]{Burchnall-Chaundy-1940}. Applying \eqref{Def-SymbolicO-1a} to \eqref{Prop-DecomF-Remark-1} and using \cite[p. 251, Eq. (17)]{Burchnall-Chaundy-1940}
	\begin{equation}\label{Prop-DecomF-Remark-2}
	(-\delta_i)_n~
	{}_{2}F_{1}\left[\begin{matrix}
		A, B\\
		C
	\end{matrix};z_i\right]
	=(-1)^n\frac{(A)_n(B)_n}{(C)_n}z_i^n ~{}_{2}F_{1}\left[\begin{matrix}
		A+n, B+n\\
		C+n
	\end{matrix};z_i\right],
	\end{equation}
	we obtain the famous result \cite[p. 252]{Burchnall-Chaundy-1940}
	\begin{align*}
		F_K^{(2)}\left[\begin{matrix}
			a_1, b_1, a_2 \\ 
			c_1, c_2
		\end{matrix}; z_1,z_2\right]
		&=\sum_{n=0}^{\infty}\frac{(a_1)_n(a_2)_n(b_1)_n}{(c_1)_n(c_2)_n}\frac{(z_1z_2)^n}{n!}\\
		&\hspace{0.5cm}\cdot {}_{2}F_{1}\left[\begin{matrix}
			a_1+n, b_1+n\\
			c_1+n
		\end{matrix};z_1\right]
		{}_{2}F_{1}\left[\begin{matrix}
		a_2+n, b_1+n\\
		c_2+n
		\end{matrix};z_2\right].
	\end{align*}
	However, it should be noted that when $L=3$, a similar method cannot yield a concise result, because we would have to compute:
	\[
	(-\delta_i)_n~z_i^m
	{}_{2}F_{1}\left[\begin{matrix}
		A, B\\
		C
	\end{matrix};z_i\right],
	\]
	and the extra factor $z_i^n$ prevents it from having a simple expression like \eqref{Prop-DecomF-Remark-2}.	
\end{remark}

Next, let us consider the following symbolic operators introduced by Hasanov \emph{et al.} \cite[pp. 2--3]{Hasanov-et al-2010}
\begin{equation}\label{Def-SymbolicO-3}
	H_{z_j}(\alpha,\beta):= \frac{\Gamma(\alpha+\delta_j)\Gamma(\beta)}{\Gamma(\alpha)\Gamma(\beta+\delta_j)}
\end{equation}
and
\begin{equation}\label{Def-SymbolicO-4}
	\bar{H}_{z_j}(\alpha,\beta):= \frac{\Gamma(\alpha)\Gamma(\beta+\delta_j)}{\Gamma(\beta)\Gamma(\alpha+\delta_j)},
\end{equation}
where $\delta_j=z_j\dfrac{\partial}{\partial z_j}$. Using \eqref{Th-GaussianST}, we can also write
\begin{equation}\label{Def-SymbolicO-3a}
	H_{z_j}(\alpha,\beta)=\sum_{n=0}^{\infty}
	\frac{(\beta-\alpha)_n(-\delta_j)_n}{(\beta)_n n!}
\end{equation}
and
\[
	\bar{H}_{z_j}(\alpha,\beta)=\sum_{n=0}^{\infty} 
	\frac{(\beta-\alpha)_n(-\delta_j)_n}{(1-\alpha-\delta_j)_n n!}.
\]

In \cite{Hasanov-et al-2010}, $9$ operator identities for the ${}_{4}F_{3}$ function are established, and from these, $11$ decomposition formulas related to ${}_{4}F_{3}$ are derived. For further studies and applications of operators \eqref{Def-SymbolicO-3} and \eqref{Def-SymbolicO-4}, interested reader may refer to \cite{Choi-Hasanov-2011}, \cite{Choi-Hasanov-2012}, \cite{Ergashev et al-2024}, \cite{Hasanov-Ergashev-2023}, and \cite{Komilova-et al-2024}.

In what follows, we apply the operator $H_{z_j}(\alpha,\beta)$ to the function $F_K^{(L)}$ and obtain the following new decomposition formula for $F_K^{(L)}$.

\begin{theorem}\label{Th-DecomF}
	Let $1\leq r\leq L-1$ and $\mathbf{z}\in\mathbb{D}_L$. We have
	\begin{align}\label{Th-DecomF-1}
		&F_K^{(L)}\left[
			\begin{matrix}
				a_1, b_1,\cdots,b_{L-1}, a_2 \\ 
				c_1,\cdots,c_L
			\end{matrix}; \mathbf{z}
			\right]
		=\sum_{n_1,\cdots,n_r=0}^{\infty} 
		(a_1)_{n_1}(b_r)_{n_r}
		\prod_{j=1}^{r} 
		\frac{(c_j-d_j)_{n_j}(-z_j)^{n_j}}{(c_j)_{n_j}(d_j)_{n_j} n_j!}\cdot \prod_{i=2}^{r}(b_{i-1})_{n_{i-1}+n_i}\notag\\
		&\hspace{0.5cm}\cdot F_K^{(L)}\left[
		\begin{matrix}
			a_1+n_1, b_1+n_1+n_2,\cdots,b_{r-1}+n_{r-1}+n_r, b_{r}+n_r, b_{r+1}, \cdots, b_{L-1}, a_2 \\ 
			d_1+n_1,\cdots,d_r+n_r, c_{r+1},\cdots,c_L
		\end{matrix}; \mathbf{z}
		\right].
	\end{align}
\end{theorem}
\begin{proof}
	By \eqref{Def-SymbolicO-3a}, we have
	\[
	H_{z_j}(d_j,c_j)
	=\sum_{n_j=0}^{\infty}\frac{(c_j-d_j)_{n_j} (-\delta_j)_{n_j}}{(c_j)_{n_j} n_j!},
	\]
	and thus
	\begin{equation}\label{Th-DecomF-Proof-1}
	\prod_{j=1}^{r}H_{z_j}(d_j,c_j)
	=\sum_{n_1,\cdots,n_r=0}^{\infty}
	\prod_{j=1}^{r} \frac{(c_j- d_j)_{n_j}(-\delta_j)_{n_j}}{(c_j)_{n_j} {n_j}!}.
	\end{equation}
	
	Let us first establish the following relation: 
	\begin{equation}\label{Th-DecomF-Proof-2}
		F_K^{(L)}\left[
		\begin{matrix}
			a_1, b_1,\cdots,b_{L-1}, a_2 \\ 
			c_1,\cdots,c_L
		\end{matrix}; \mathbf{z}
		\right]=H_{z_1}(d_1,c_1)\cdots H_{z_r}(d_r,c_r) 
		F_K^{(L)}\left[
		\begin{matrix}
			a_1, b_1,\cdots,b_{L-1}, a_2 \\ 
			d_1,\cdots,d_r,c_{r+1},\cdots,c_{L}
		\end{matrix}; \mathbf{z}
		\right].
	\end{equation}
	Actually, we have
	\begin{align*}
	H_{z_r}(d_r,c_r)
	F_K^{(L)}\left[
	\begin{matrix}
		a_1, b_1,\cdots,b_{L-1}, a_2 \\ 
		d_1,\cdots,d_r,c_{r+1},\cdots,c_{L}
	\end{matrix}; \mathbf{z}
	\right]
	&=\sum_{\mathbf{n}=\mathbf{0}}^{\infty} \frac{(a_1)_{n_1}(b_1)_{n_1+n_2}\cdots(b_{L-1})_{n_{L-1}+n_L}(a_2)_{n_L}}{(d_1)_{n_1}\cdots(d_r)_{n_r}(c_{r+1})_{n_{r+1}}\cdots(c_L)_{n_L}}\\
	&\cdot \frac{z_1^{n_1}\cdots z_{r-1}^{n_{r-1}}}{n_1!\cdots n_{r-1}!}\cdot\frac{1}{n_r!}\cdot\frac{z_{r+1}^{n_{r+1}}\cdots z_L^{n_L}}{n_{r+1}!\cdots n_L!}
	\cdot H_{z_r}(d_r,c_r)\left\{z_r^{n_r}\right\}.
	\end{align*}
	Using \eqref{Def-SymbolicO-3} and \eqref{PooleF-1} gives immediately 
	\[
	 H_{z_r}(d_r,c_r)\left\{z_r^{n_r}\right\}
	 =\frac{(d_r)_{n_r}}{(c_r)_{n_r}}z_r^{n_r}.
	\]
	Thus, we obtain 
	\begin{align*}
		H_{z_r}(d_r,c_r)
		F_K^{(L)}\left[
		\begin{matrix}
			a_1, b_1,\cdots,b_{L-1}, a_2 \\ 
			d_1,\cdots,d_r,c_{r+1},\cdots,c_{L}
		\end{matrix}; \mathbf{z}
		\right]
		&=F_K^{(L)}\left[
		\begin{matrix}
			a_1, b_1,\cdots,b_{L-1}, a_2 \\ 
			d_1,\cdots,d_{r-1},c_r,\cdots,c_{L}
		\end{matrix}; \mathbf{z}
		\right].
	\end{align*}
	Continuing this process, we get \eqref{Th-DecomF-Proof-2}.

	Based on Eq. \eqref{Th-DecomF-Proof-2}, it follows from \eqref{Th-DecomF-Proof-1} that
	\[
	F_K^{(L)}\left[
	\begin{matrix}
		a_1, b_1,\cdots,b_{L-1}, a_2 \\ 
		c_1,\cdots,c_L
	\end{matrix}; \mathbf{z}
	\right]
	=\sum_{n_1,\cdots,n_r=0}^{\infty}\prod_{j=1}^{r} 
	\frac{(c_j-d_j)_{n_j}(-\delta_j)_{n_j}}{(c_j)_{n_j} n_j!}
	F_K^{(L)}\left[
	\begin{matrix}
		a_1, b_1,\cdots,b_{L-1}, a_2 \\ 
		d_1,\cdots,d_r, c_{r+1},\cdots,c_L
	\end{matrix}; \mathbf{z}
	\right].
	\]
	Using \eqref{PooleF-1}, we obtain
	\begin{align*}
	F_K^{(L)}\left[
	\begin{matrix}
		a_1, b_1,\cdots,b_{L-1}, a_2 \\ 
		c_1,\cdots,c_L
	\end{matrix}; \mathbf{z}
	\right]
	&=\sum_{n_1,\cdots,n_r=0}^{\infty}\prod_{j=1}^{r} 
	\frac{(c_j-d_j)_{n_j}(-z_j)^{n_j}}{(c_j)_{n_j} n_j!}\notag\\
	&\hspace{0.5cm}\cdot 
	\frac{\partial^{n_1 + \cdots + n_r}}{\partial z_1^{n_1} \cdots \partial z_r^{n_r}}
	F_K^{(L)}\left[
	\begin{matrix}
		a_1, b_1,\cdots,b_{L-1}, a_2 \\ 
		d_1,\cdots,d_r, c_{r+1},\cdots,c_L
	\end{matrix}; \mathbf{z}
	\right],
	\end{align*}
	where
	\begin{align*}
		&\frac{\partial^{n_1 + \cdots + n_r}}{\partial z_1^{n_1} \cdots \partial z_r^{n_r}}
		F_K^{(L)}\left[
		\begin{matrix}
			a_1, b_1,\cdots,b_{L-1}, a_2 \\ 
			d_1,\cdots,d_r, c_{r+1},\cdots,c_L
		\end{matrix}; \mathbf{z}
		\right]=(a_1)_{n_1}(b_r)_{n_r}
		\prod_{j=1}^{r-1}(b_j)_{n_j+n_{j+1}}
		\\
		&\hspace{0.5cm}\cdot \prod_{j=1}^{r}\frac{1}{(d_j)_{n_j}}\cdot F_K^{(L)}\left[
		\begin{matrix}
		a_1+n_1, b_1+n_1+n_2,\cdots,b_{r-1}+n_{r-1}+n_r, b_{r}+n_r, b_{r+1}, \cdots, b_{L-1}, a_2 \\ 
		d_1+n_1,\cdots,d_r+n_r, c_{r+1},\cdots,c_L
		\end{matrix}; \mathbf{z}
		\right].
	\end{align*}
	This completes the proof. 
\end{proof}

\emph{An alternative proof of Theorem \ref{Th-DecomF}.} For convenience, let us denote the right-hand side of \eqref{Th-DecomF-1} by $S$. With the help of \eqref{IntegralR-FKL}, we can express $S$ as
\begin{align}\label{Th-DecomF-Proof2-1}
	S&=\sum_{n_1,\cdots,n_r=0}^{\infty}
		(a_1)_{n_1}(b_r)_{n_r}\prod_{j=1}^{r} 
		\frac{(c_j-d_j)_{n_j}(-z_j)^{n_j}}{(c_j)_{n_j}(d_j)_{n_j}n_j!}
		\prod_{i=2}^{r}(b_{i-1})_{n_{i-1}+n_i}\notag\\
	&\hspace{0.5cm}\cdot\int_{0}^{1}\int_{0}^{1}
		(1-t_1z_1)^{-b_1-n_1-n_2}
		(1-t_Lz_L)^{-b_{L-1}}\notag\\
	&\hspace{0.5cm}\cdot 
	F_K^{(L-2)}\Bigg[\begin{matrix} 
			b_1+n_1+n_2,\textcolor{orange}{\cdots,} b_{r-1}+n_{r-1}+n_r, b_r+n_r,b_{r+1},\cdots, b_{L-1} \\ d_2+n_2,\cdots,d_r+n_r,c_{r+1},\cdots, c_{L-1} \end{matrix};\notag\\
	&\hspace{1.7cm}\frac{z_2}{1-t_1z_1},z_3,\cdots,z_{L-2},\frac{z_{L-1}}{1-t_Lz_L}\Bigg]
	\mathrm{d}\mu_{a_1+n,d_1-a_1}(t_1)\mathrm{d}\mu_{a_2,c_L-a_2}(t_L).
\end{align}
Expressing the $F_{K}^{(L-2)}$ function in \eqref{Th-DecomF-Proof2-1} as a series and interchanging the order of integration and summation, we have
\begin{align}\label{Th-DecomF-Proof2-2}
	S&=\sum_{n_1,\cdots,n_r=0}^{\infty}
	(a_1)_{n_1}(b_r)_{n_r}\prod_{j=1}^{r}\frac{(c_j-d_j)_{n_j}(-z_j)^{n_j}}{(c_j)_{n_j}(d_j)_{n_j}n_j!}
	\prod_{i=2}^{r}(b_{i-1})_{n_{i-1}+n_i}\notag\\
	&\hspace{0.5cm}\cdot\sum_{m_2,\cdots,m_{L-1}=0}^{\infty}
	\frac{(b_1+n_1+n_2)_{m_2}
	(b_2+n_2+n_3)_{m_2+m_3}
	\cdots (b_{r-1}+n_{r-1}+n_r)_{m_{r-1}+m_r}
	}{(d_2+n_2)_{m_2}\cdots(d_r+n_r)_{m_r}}\notag\\
	&\hspace{0.5cm}\cdot \frac{(b_r+n_r)_{m_r+m_{r+1}}
		(b_{r+1})_{m_{r+1}+m_{r+2}}
		\cdots(b_{L-2})_{m_{L-2}+m_{L-1}}
		(b_{L-1})_{m_{L-1}}}{(c_{r+1})_{m_{r+1}}\cdots (c_{L-1})_{m_{L-1}} m_2!\cdots m_{L-1}!}
		z_2^{m_2}\cdots z_{L-1}^{m_{L-1}}\notag\\
	&\hspace{0.5cm}\cdot \int_{0}^{1}(1-t_1z_1)^{-b_1-n_1-n_2-m_2}\mathrm{d}\mu_{a_1+n_1,d_1-a_1}(t_1)\int_{0}^{1}(1-t_Lz_L)^{-b_{L-1}-m_{L-1}}\mathrm{d}\mu_{a_2,c_L-a_2}(t_L).
\end{align}
	
Applying \eqref{IntegralR-FKL-L1} to \eqref{Th-DecomF-Proof2-2} gives
	\begin{align*}
	S&=\sum_{n_1,\cdots,n_r=0}^{\infty}
		\prod_{j=1}^{r} 
		\frac{(c_j-d_j)_{n_j}(-z_j)^{n_j}}{(c_j)_{n_j}n_j!}
		\sum_{m_2,\cdots,m_{L-1}=0}^{\infty}
		\frac{(a_1)_{n_1}}{(d_1)_{n_1}}\\
		&\hspace{0.5cm}\cdot
		{}_{2}F_1\left[\begin{matrix} 
			b_1+n_1+n_2+m_2, a_1+n_1 \\ 
			d_1+n_1
			\end{matrix}; z_1 \right]
		{}_{2}F_1\left[\begin{matrix}
			b_{L-1}+m_{L-1}, a_2 \\ 
			c_L 
			\end{matrix}; z_L \right]\\
		&\hspace{0.5cm}\cdot\frac{(b_1)_{n_1+n_2+m_2}(b_2)_{n_2+n_3+m_2+m_3}\cdots(b_{r-1})_{n_{r-1}+n_r+m_{r-1}+m_r}(b_r)_{n_r+m_r+m_{r+1}}}{(d_2)_{n_2+m_2}\cdots(d_r)_{n_r+m_r}}\\
		&\hspace{0.5cm}\cdot
		\frac{(b_{r+1})_{m_{r+1}+m_{r+2}}\cdots(b_{L-2})_{m_{L-2}+m_{L-1}}(b_{L-1})_{m_{L-1}}}{(c_{r+1})_{m_{r+1}}\cdots(c_{L-1})_{m_{L-1}}m_2!\cdots m_{L-1}!}
		z_2^{m_2}\cdots z_{L-1}^{m_{L-1}}.
	\end{align*}
	Expanding the two ${}_{2}F_{1}$ functions into series and simplifying the resulting expression, we obtain
	\begin{align*}
		S&=\sum_{n_1,\cdots,n_r=0}^{\infty}
		\prod_{j=1}^{r} 
		\frac{(c_j-d_j)_{n_j}(-z_j)^{n_j}}{(c_j)_{n_j}n_j!}\sum_{m_1,\cdots,m_L=0}^{\infty}
		(a_1)_{n_1+m_1}(a_2)_{m_L}\\
		&\hspace{0.5cm}\cdot\frac{(b_1)_{n_1+n_2+m_1+m_2}\cdots(b_{r-1})_{n_{r-1}+n_r+m_{r-1}+m_r}(b_r)_{n_r+m_r+m_{r+1}}}{(d_1)_{n_1+m_1}\cdots(d_r)_{n_r+m_r}}\\
		&\hspace{0.5cm}\cdot
		\frac{(b_{r+1})_{m_{r+1}+m_{r+2}}\cdots(b_{L-2})_{m_{L-2}+m_{L-1}}(b_{L-1})_{m_{L-1}+m_L}}{(c_{r+1})_{m_{r+1}}\cdots(c_L)_{m_L}m_1!\cdots m_L!}
		z_1^{m_1}\cdots z_L^{m_L}.
	\end{align*}
	
	Setting $n_j+m_j=\ell_j$ $(n_j=\ell_j-m_j)$ for $j=1,\cdots, r$ and interchanging the order of summation, we obtain
	\begin{align*}
		S&=\sum_{\ell_1,\cdots,\ell_r=0}^{\infty}
		\sum_{m_1=0}^{\ell_1}\cdots\sum_{m_r=0}^{\ell_r}
		\prod_{j=1}^{r} 
		\frac{(c_j-d_j)_{\ell_j-m_j}(-1)^{\ell_j-m_j}z_j^{\ell_j}}{(c_j)_{\ell_j-m_j}m_j!(\ell_j-m_j)!}\sum_{m_{r+1},\cdots,m_L=0}^{\infty}
		(a_1)_{\ell_1}(a_2)_{m_L}\\
		&\hspace{0.5cm}\cdot\frac{(b_1)_{\ell_1+\ell_2}\cdots(b_{r-1})_{\ell_{r-1}+\ell_r}(b_r)_{\ell_r+m_{r+1}}}{(d_1)_{\ell_1}\cdots(d_r)_{\ell_r}}\\
		&\hspace{0.5cm}\cdot
		\frac{(b_{r+1})_{m_{r+1}+m_{r+2}}\cdots(b_{L-2})_{m_{L-2}+m_{L-1}}(b_{L-1})_{m_{L-1}+m_L}}{(c_{r+1})_{m_{r+1}}\cdots(c_L)_{m_L}m_{r+1}!\cdots m_L!}
		z_{r+1}^{m_{r+1}}\cdots z_L^{m_L}.
	\end{align*}
	Applying the identity $(-\ell)_{m}=\dfrac{(-1)^{m}\ell!}{(\ell-m)!}$ together with \eqref{Th-ChuVandermonde} yields
	\begin{align*}
	\sum_{m_j=0}^{\ell_j}\frac{(c_j-d_j)_{\ell_j-m_j}(-1)^{\ell_j-m_j}}{(c_j)_{\ell_j-m_j}m_j!(\ell_j-m_j)!}
	=\frac{1}{\ell_j!}\frac{(d_j)_{\ell_j}}{(c_j)_{\ell_j}}.
	\end{align*}
	Thus we have
	\begin{align*}
		G&=\sum_{\ell_1,\cdots,\ell_r=0}^{\infty}
		\sum_{m_{r+1},\cdots,m_L=0}^{\infty}
		(a_1)_{\ell_1}\frac{(b_1)_{\ell_1+\ell_2}\cdots(b_{r-1})_{\ell_{r-1}+\ell_r}(b_r)_{\ell_r+m_{r+1}}}{(c_1)_{\ell_1}\cdots(c_r)_{\ell_r}}\\
		&\hspace{0.5cm}\cdot
		\frac{(b_{r+1})_{m_{r+1}+m_{r+2}}\cdots(b_{L-2})_{m_{L-2}+m_{L-1}}(b_{L-1})_{m_{L-1}+m_L}}{(c_{r+1})_{m_{r+1}}\cdots(c_L)_{m_L}
			\ell_1!\cdots \ell_r! m_{r+1}!\cdots m_L!}(a_2)_{m_L}
		z_1^{\ell_1}\cdots z_r^{\ell_r}
		z_{r+1}^{m_{r+1}}\cdots z_L^{m_L}\\
	&=F_K^{(L)}\left[
	\begin{matrix}
		a_1, b_1,\cdots,b_{L-1}, a_2 \\ 
		c_1,\cdots,c_L
	\end{matrix}; \mathbf{z}
	\right].
	\end{align*}
	This completes the proof. 
	
\begin{remark}
	If we take $L=2$ and $r=1$ in \eqref{Th-DecomF-1}, we obtain the following for Appell function $F_2$: 
	\[
		F_K^{(2)}\left[
		\begin{matrix}
			a_1, b_1, a_2 \\ 
			c_1, c_2
		\end{matrix}; z_1,z_2
		\right]
		=\sum_{n=0}^{\infty} 
		\frac{(a_1)_n (b_1)_n (c_1-d_1)_n}{(c_1)_n(d_1)_n}
		\frac{(-z_1)^n}{n!}
		 F_K^{(2)}\left[
		\begin{matrix}
			a_1+n, b_1, a_2 \\ 
			d_1+n, c_2
		\end{matrix}; z_1,z_2
		\right].
	\]
\end{remark}

\section{Concluding remarks}\label{Sect-5}

In this paper, we study the multivariate hypergeometric function $F_K^{(L)}$, defined by \eqref{Def-MultivariateFK}, as a new special function for the first time, and establish some of its fundamental properties. In particular, we present different approaches to establishing Proposition \ref{Prop-ReductionF-ConvergenceC} and Theorem \ref{Th-DecomF}, offering alternative perspectives that deepen our understanding of the function. However, many problems concerning $F_K^{(L)}$ remain to be studied. Here, we propose two directions worthy of future research.
\begin{itemize}
	\item[(1)] For hypergeometric functions, various types of integral representations offer different advantages in addressing different problems. In this paper, we have focused exclusively on Euler-type integral representations. Other types of integral representations, such as Laplace-type and Mellin--Barnes-type representations, therefore warrant further investigation. 
	\item[(2)] The asymptotic behavior of hypergeometric functions has long been an important topic in the theory of special functions. For the case $L=3$, a considerable number of results have recently been obtained. Considering the practical applications of the $F_K^{(L)}$ function in physics and the inherent complexity of multidimensional problems, studying its various types of asymptotic behavior will definitely broaden our horizons in both theory and application.
\end{itemize}

Finally, we would like to mention the relationship between $F_K^{(L)}$ and another new multivariate hypergeometric function. In 2023, Alkalaev and Mandraygin introduce the following multivariate hypergeometric function \cite[p. 2, Eq. (1.2)]{Alkalaev-Mandrygin-2023}: 
\[
	F_N^{(L)}\left[\begin{matrix}
	\mathbf{b}\\
	\mathbf{c}
	\end{matrix}; \mathbf{z}\right]
	:=\sum_{\mathbf{n}=\mathbf{0}}^{\infty}
	\frac{(b_1)_{n_1+n_2}(b_2)_{n_2+n_3}\cdots
		(b_{L-1})_{n_{L-1}+n_L} (b_{L})_{n_{L}+n_{1}}}{(c_1)_{n_1}\cdots(c_L)_{n_L}}
	\frac{\mathbf{z}^{\mathbf{n}}}{\mathbf{n}!}.
\]
They call $F_N^{(L)}$ the \emph{necklace function} because it describes the bare global conformal block in the necklace channel. When $L=1$, it may reduce to the Gauss hypergeometric function ${}_{2}F_{1}$ of the form
\[
{}_{2}F_{1}\left[\begin{matrix}
	\frac{a_1}{2},\frac{1+a_1}{2}\\
	c_1\end{matrix};4z_1\right]. 
\]
When $L=2$, it becomes Appell's function $F_4$ \cite[p. 23, Eq. (5)]{Srivastava-Karlsson-Book-1985}. When $L=3$, it becomes Srivastava's $H_B$ function \cite[p. 43, Eq. (12)]{Srivastava-Karlsson-Book-1985}. 
Here, we likewise suggest carrying out a systematic study of the 
$F_N^{(L)}$ function. Actually, in view of \eqref{DirichletM-1}, it is not difficult to see that
\[
\int_0^1 t^{n_1}(1-t)^{n_L}\mathrm{d}\mu_{a_1,b_L-a_1}(t)
=\frac{(a_1)_{n_1}(b_L-a_1)_{n_L}}{(b_L)_{n_1+n_L}},
\]
where $\Re(b_L)>\Re(a_1)>0$ and $n_1,n_L\in\mathbb{Z}_{\geq0}$. Then we have 
\[
\int_0^1
F_N^{(L)}\left[\begin{matrix}
	\mathbf{b}\\
	\mathbf{c}
\end{matrix}; tz_1,z_2,\cdots,z_{L-1},(1-t)z_L\right]\mathrm{d}\mu_{a_1,b_L-a_1}(t)
=F_K^{(L)}\left[\begin{matrix}
	a_1, b_1, \cdots, b_{L-1}, b_L-a_1\\ 
	c_1, \cdots, c_L
\end{matrix}; \mathbf{z}\right].
\]
From this, it can be seen that $F_N^{(L)}$ and $F_K^{(L)}$ are very closely related functions, and perhaps they can be studied within the same framework.

\appendix
\section{The use of Horn's theorem}\label{Appendix-1}

We adopt the following version of the Horn-type theorem given by Sadykov and Tsikh (see \cite[p. 148, Theorem 4.1.6]{Sadykov-Tsikh-2014}, see also \cite[p. 223, Theorem 1]{Semusheva-Tsikh-2009}).

\begin{theorem}
	Consider the multiple power series of $L$ variables 
	\[
	H(\mathbf{z})=\sum_{\mathbf{n}=\mathbf{0}}^{\infty}f(\mathbf{n})\mathbf{z}^{\mathbf{n}},
	\]
	where
	\[
	f(\mathbf{n}):=\frac{\displaystyle\prod_{i=1}^{p}\Gamma(\langle\mathbf{A}_i,\mathbf{n}\rangle+a_i)}{\displaystyle\prod_{j=1}^{q}\Gamma(\langle\mathbf{B}_i,\mathbf{n}\rangle+b_j)}.
	\]
	Here, $a_i,b_j\in\mathbb{C}$, $\mathbf{A}_i=(A_{i,1},\cdots,A_{i,L})\in \mathbb{R}^L$, and $\mathbf{B}_i=(B_{i,1},\cdots,B_{i,L})\in \mathbb{R}^L$. 
	Define 
	\[
	F_i(\mathbf{q}):=\lim_{t\rightarrow+\infty}\frac{f(t\mathbf{q}+\mathbf{e}_i)}{f(t\mathbf{q})}, \quad i=1,\cdots,L,
	\]
	where $\mathbf{q}=(q_1,\cdots,q_L)$.
	
	For any nonempty subset $J\subseteq I=\{1,\cdots,L\}$ of cardinality $\mathrm{card}(J)$, define vector-valued function
	\[
	F_J(\mathbf{q}_J,\mathbf{0}_{I\setminus J})
	=(F_j(\mathbf{q}_J,\mathbf{0}_{I\setminus J}))_{j\in J}:
	\mathbb{R}_{>0}^{\mathrm{card}(J)}\mapsto\mathbb{R}_{>0}^{\mathrm{card}(J)},
	\]
	where $(\mathbf{q}_J,\mathbf{0}_{I\setminus J})$ is the $L$-dimensional vector whose coordinates with $j\in J$ are $q_j$, and whose remaining coordinates are zero. Then the region of convergence of $H(\mathbf{z})$ is given by
	\[
	G=\bigcap_{\emptyset\neq J\subseteq I}G_J
	=\bigcap_{k=1}^{L}\left(\bigcap_{\substack{J\subseteq I\\\mathrm{card}(J)=k}}G_J\right), 
	\]
	where
	\[
	G_J:=\left\{\mathbf{z}\in\mathbb{C}^L;~\forall \mathbf{q}_J\in\mathbb{R}_{> 0}^{\mathrm{card}(J)}~\exists j\in J:~ |z_j|<\frac{1}{|F_j(\mathbf{q}_J,\mathbf{0}_{I\setminus J})|}\right\}.
	\]
\end{theorem}

When $L=4$, we need to consider the intersection of $2^4-1=15$ sets. Furthermore, a complete statement of the Horn-type theorem for the 4-dimensional case can also be found in \cite[pp. 3--4]{Comeau et al-2020}.

Now, let us consider 
\[
	F_K^{(4)}\left[\begin{matrix}
		a_1, b_1, b_2, b_3, a_2 \\ 
		c_1, c_2, c_3, c_4
	\end{matrix}; z_1,z_2,z_3,z_4\right]
	=\sum_{n_1, \cdots, n_4=0}^{\infty}
	f(n_1,n_2,n_3,n_4) z_1^{n_1}z_2^{n_2}z_3^{n_3}z_4^{n_4},
\]
where
\[
f(n_1,n_2,n_3,n_4):=
\frac{(a_1)_{n_1}(b_1)_{n_1+n_2}(b_2)_{n_2+n_3}(b_3)_{n_3+n_4}(a_2)_{n_4}}{(c_1)_{n_1}(c_2)_{n_2}(c_3)_{n_3}(c_4)_{n_4}n_1!n_2!n_3!n_4!}.
\]
It is not difficult to check that
\begin{align*}
	F_1(q_1,q_2,q_3,q_4)&=\frac{q_1+q_2}{q_1},\\
	F_2(q_1,q_2,q_3,q_4)&=\frac{(q_1+q_2)(q_2+q_3)}{q_2^2},\\
	F_3(q_1,q_2,q_3,q_4)&=\frac{(q_2+q_3)(q_3+q_4)}{q_3^2},\\
	F_4(q_1,q_2,q_3,q_4)&=\frac{q_3+q_4}{q_4}.
\end{align*}
Let $r_j=|z_j|$ $(j=1,\cdots,4)$. If $J=\{1\}$, then 
\[
G_{\{1\}}=\left\{\mathbf{z}\in\mathbb{C}^4;~\forall q_1\in\mathbb{R}_{> 0}:~ |z_1|<\frac{1}{|F_1(q_1,0,0,0)|}\right\}
=\left\{\mathbf{z}\in\mathbb{C}^4;r_1<1\right\}.
\]
If $J=\{2,3\}$, then 
\begin{align*}
G_{\{2,3\}}&=\left\{\mathbf{z}\in\mathbb{C}^4;~\forall (q_2,q_3)\in\mathbb{R}_{>0}^{2}~\exists j\in \{2,3\}:~ |z_j|<\frac{1}{|F_j(0,q_2,q_3,0)|}\right\}\\
&=\left\{\mathbf{z}\in\mathbb{C}^4;~\forall (q_2,q_3)\in\mathbb{R}_{> 0}^{2}:|z_2|<\frac{1}{|F_2(0,q_2,q_3,0)|}\vee|z_3|<\frac{1}{|F_3(0,q_2,q_3,0)|}\right\}\\
&=\left\{\mathbf{z}\in\mathbb{C}^4; r_2<\frac{q_2}{q_2+q_3}\vee r_3<\frac{q_3}{q_2+q_3}\right\}.
\end{align*}
Proceeding in this manner, we obtain the full representation of these 15 sets:
\begin{align*}
	G_{\{1\}}&=\left\{\mathbf{z}\in\mathbb{C}^4;r_1<1\right\},\\
	G_{\{2\}}&=\left\{\mathbf{z}\in\mathbb{C}^4;r_2<1\right\},\\
	G_{\{3\}}&=\left\{\mathbf{z}\in\mathbb{C}^4;r_3<1\right\},\\
	G_{\{4\}}&=\left\{\mathbf{z}\in\mathbb{C}^4;r_4<1\right\},\\
	G_{\{1,2\}}&=\left\{\mathbf{z}\in\mathbb{C}^4;r_1<\frac{q_1}{q_1+q_2}\vee r_2<\frac{q_2}{q_1+q_2}\right\},\\
	G_{\{1,3\}}&=\left\{\mathbf{z}\in\mathbb{C}^4;r_1<1\vee r_3<1\right\},\\
	G_{\{1,4\}}&=\left\{\mathbf{z}\in\mathbb{C}^4;r_1<1\vee r_4<1\right\},\\
	G_{\{2,3\}}&=\left\{\mathbf{z}\in\mathbb{C}^4;r_2<\frac{q_2}{q_2+q_3}\vee r_3<\frac{q_3}{q_2+q_3}\right\},\\
	G_{\{2,4\}}&=\left\{\mathbf{z}\in\mathbb{C}^4;r_2<1\vee r_4<1\right\},\\
	G_{\{3,4\}}&=\left\{\mathbf{z}\in\mathbb{C}^4;r_3<\frac{q_3}{q_3+q_4}\vee r_4<\frac{q_4}{q_3+q_4}\right\},\\
	G_{\{1,2,3\}}&=\left\{\mathbf{z}\in\mathbb{C}^4;r_1<\frac{q_1}{q_1+q_2}\vee r_2<\frac{q_2^2}{(q_1+q_2)(q_2+q_3)} \vee r_3<\frac{q_3}{q_2+q_3}\right\},\\
	G_{\{1,2,4\}}&=\left\{\mathbf{z}\in\mathbb{C}^4;r_1<\frac{q_1}{q_1+q_2}\vee r_2<\frac{q_2}{q_1+q_2} \vee r_4<1\right\},\\
	G_{\{1,3,4\}}&=\left\{\mathbf{z}\in\mathbb{C}^4;r_1<1\vee r_3<\frac{q_3}{q_3+q_4} \vee r_4<\frac{q_4}{q_3+q_4}\right\},\\
	G_{\{2,3,4\}}&=\left\{\mathbf{z}\in\mathbb{C}^4;r_2<\frac{q_2}{q_2+q_3}\vee r_3<\frac{q_3^2}{(q_2+q_3)(q_3+q_4)} \vee r_4<\frac{q_4}{q_3+q_4}\right\},\\
	G_{\{1,2,3,4\}}&=\left\{\mathbf{z}\in\mathbb{C}^4;r_1<\frac{q_1}{q_1+q_2}\vee r_2<\frac{q_2^2}{(q_1+q_2)(q_2+q_3)} \vee r_3<\frac{q_3^2}{(q_2+q_3)(q_3+q_4)}\vee r_4<\frac{q_4}{q_3+q_4}\right\}.
\end{align*}
Thus, the region of convergence of $F_K^{(4)}$ is the intersection set of above $15$ sets. To obtain the simplified expression, we will only focus on the active inequality constraints in the sequel. Actually, from $G_{\{1\}}, \cdots, G_{\{4\}}$, we have 
\[
r_1<1, \quad r_2<1, \quad r_3<1, \quad r_4<1.
\]
From $G_{\{1,2\}}, G_{\{2,3\}}, G_{\{3,4\}}$, we obtain
\[
r_1+r_2<1, \quad r_2+r_3<1, \quad r_3+r_4<1.
\]
Note that the sets $G_{\{1,3\}}, G_{\{1,4\}}, G_{\{2,4\}}$ does not provide active inequality constraints. From $G_{\{1,2,3\}}, G_{\{2,3,4\}}$, we have
\[
r_2<(1-r_1)(1-r_3), \quad r_3<(1-r_2)(1-r_4). 
\]
The sets $G_{\{1,2,4\}}, G_{\{1,3,4\}}$ does not provide active inequality constraints. From $G_{\{1,2,3,4\}}$, we obtain
\[
r_2r_3<(1-r_1-r_2)(1-r_3-r_4),
\]
that is, 
\begin{equation}\label{Appendix-1-1}
1-r_1-r_2-r_3-r_4+r_1r_3+r_1r_4+r_2r_4>0.
\end{equation}
Under the condition $r_1<1$ and $r_4<1$, \eqref{Appendix-1-1} is equivalent to 
\[
\frac{r_2}{1-r_1}+\frac{r_3}{1-r_4}<1.
\]
Therefore, the result derived by Horn's method is completely consistent with the convergence domain $\mathbb{D}_4$ of $F_K^{(4)}$ that we obtained in \eqref{Cond-FK4}.

\section*{Funding}

No funding was received for conducting this study.

\section*{Contributions}

Both authors contributed equally to the preparation of this paper.

\section*{Conflicts of interest}

The authors declare that there is no conflict of interest.

\section*{Data Availability}

This manuscript has no associated data.

\section*{Acknowledgements}

We are grateful to Peng-Cheng Hang for his careful reading of the initial manuscript and providing valuable suggestions.


\begin{thebibliography}{99}
\footnotesize
\setlength{\itemsep}{-2pt}

\bibitem{Alkalaev-Mandrygin-2023}
K. Alkalaev, S. Mandrygin, 
Torus shadow formalism and exact global conformal blocks, 
J. High Energy Phys. 2023 (2023), Paper No. 157, 26 pp. 
\url{https://doi.org/10.1007/JHEP11(2023)157} 

\bibitem{Ammon et al-2025}
M. Ammon, J. Hollweck, T. H\"{o}ssel, K. W\"{o}lfl,  
Conformal blocks in two and four dimensions from oscillator representations, 
J. High Energy Phys. 2025 (2025), Paper No. 91, 34 pp.
\url{https://doi.org/10.1007/JHEP05(2025)091}

\bibitem{Antonova-Dmytryshyn-Goran-2023}
T. Antonova, R. Dmytryshyn, V. Goran, 
On the analytic continuation of Lauricella–Saran hypergeometric function $F_K(a_1,a_2,b_1,b_2;a_1,b_2,c_3;z)$, 
Mathematics 11 (2023), 4487.
\url{https://doi.org/10.3390/math11214487}

\bibitem{Belitsky-2018}
A.V. Belitsky,
Multichannel conformal blocks for scattering amplitudes,
Phys. Lett., B 780 (2018), 66--73.
\url{https://doi.org/10.1016/j.physletb.2018.02.063}

\bibitem{Bhattacharyya et al-2025}
A. Bhattacharyya, S. Ghosh, P. Nandi, S. Pal, 
3D $\mathcal{N}=1$ supergravity from Virasoro TQFT: gravitational partition function and out-of-time-order correlator, 
J. High Energy Phys. 2025 (2025), Paper No. 27, 51 pp. 
\url{https://doi.org/10.1007/JHEP02(2025)027}

\bibitem{Burchnall-Chaundy-1940}
J.L. Burchnall, T.W. Chaundy,
Expansions of Appell's double hypergeometric functions,
Q. J. Math., Oxf. Ser. 11 (1940), 249--270.
\url{https://doi.org/10.1093/qmath/os-11.1.249}

\bibitem{Burchnall-Chaundy-1941}
J.L. Burchnall, T.W. Chaundy,
Expansions of Appell's double hypergeometric functions (II),
Q. J. Math., Oxf. Ser. 12 (1941), 112--128.
\url{https://doi.org/10.1093/qmath/os-12.1.112}

\bibitem{Carlson-Book-1977}
B.C. Carlson, Special Functions of Applied Mathematics, Academic Press, New York, 1977.

\bibitem{Choi-Hasanov-2011}
J. Choi, A. Hasanov, 
Applications of the operator $H(\alpha,\beta)$ to the Humbert double hypergeometric functions, 
Comput. Math. Appl. 61 (3) (2011), 663--671.
\url{https://doi.org/10.1016/j.camwa.2010.12.012}

\bibitem{Choi-Hasanov-2012}
J. Choi, A. Hasanov, 
Certain decomposition formulas of generalized hypergeometric functions ${}_{p}F_{q}$ and some formulas of an analytic continuation of the Clausen function ${}_{3}F_{2}$, 
Commun. Korean Math. Soc. 27 (1) (2012), 107--116. 
\url{https://doi.org/10.4134/CKMS.2012.27.1.107}

\bibitem{Comeau et al-2020}
V. Comeau, J.-F. Fortin, W. Skiba, 
Further Results on a Function Relevant for Conformal Blocks, 
SIGMA 16 (2020), 124, 15 pp. 
\url{https://doi.org/10.3842/SIGMA.2020.124}

\bibitem{Dmytryshyn-Goran-2024}
R. Dmytryshyn, V. Goran, 
On the analytic extension of Lauricella-Saran's hypergeometric function $F_K$ to symmetric domains, 
Symmetry 16 (2024), 220.
\url{https://doi.org/10.3390/sym16020220}

\bibitem{Dmytryshyn-Goran-2026}
R. Dmytryshyn, V. Goran, 
On truncation error bounds of branched continued fraction expansions of some ratios of Lauricella-Saran's hypergeomtric functions $F_K$, 
Mat. Stud. 65 (1) (2026), 48--57.
\url{https://doi.org/10.30970/ms.65.1.48-57}

\bibitem{Erdelyi-1939b}
A. Erd\'{e}lyi, 
Transformation of hypergeometric integrals by means of fractional integration by parts, 
Q. J. Math., Oxf. Ser. 10 (1939), 176--189.

\bibitem{Ergashev et al-2024}
T.G. Ergashev, A. Hasanov, T.K. Yuldashev, B.G. Shamilova,
The infinite summation formulas of the multiple Lauricella's hypergeometric functions,
Trans. Natl. Acad. Sci. Azerb., Ser. Phys.-Tech. Math. Sci. 44 (2024), 22--33.
\url{https://doi.org/10.30546/2617-7900.44.4.2024.022}

\bibitem{Ferrando et al-2026}
G. Ferrando, F. Loebbert, A. Pitters, S.F. Stawinski,
Hypergeometry from $\hat{P}$-symmetry: Feynman integrals in one and two dimensions,
J. High Energy Phys. 2026 (3) (2026), Paper No. 255, 76 p.
\url{https://doi.org/10.1007/JHEP03(2026)255}

\bibitem{Guo-Luo-2026}
L.-J. Guo, M.-J. Luo, 
Erd\'{e}lyi-type integrals for $F_K$ function and their $q$-analogues, Fractal Fract. 10 (4) (2026), 225. 
\url{https://doi.org/10.3390/fractalfract10040225}

\bibitem{Hang et al-2025}
P.-C. Hang, M. Henkel, M.-J. Luo,  
Asymptotics of the Humbert functions $\Psi_1$ and $\Psi_2$, 
J. Approx. Theory 314 (2) (2026), Article ID 106233, 25 pp.  
\url{https://doi.org/10.1016/j.jat.2025.106233}

\bibitem{Hang et al-2026a}
P.-C. Hang, L. Hu, M.-J. Luo, 
Complete asymptotic expansions of the Humbert function $\Psi_1$ for two large arguments, Math. Slovaca, 76 (5) (2026), 1213--1225.
\url{https://doi.org/10.1515/ms-2026-0111}

\bibitem{Hang et al-2026b}
P.-C. Hang, L. Hu, M.-J. Luo, 
Asymptotic expansions of the Humbert function $\Phi_1$ and their applications, 2026. 
\url{https://arxiv.org/abs/2504.09280v3}

\bibitem{Hang-Luo-2025}
P.-C. Hang, M.-J. Luo,
Asymptotics of Saran's hypergeometric function $F_K$, 
J. Math. Anal. Appl. 541 (2) (2025), Article ID 128707, 19 pp.
\url{https://doi.org/10.1016/j.jmaa.2024.128707}

\bibitem{Hang-Luo-2024}
P.-C. Hang, M.-J. Luo, 
Asymptotics of the Humbert function $\Psi_1$ for two large arguments, contribution to \emph{Special Issue on Asymptotics and Applications of Special Functions in Memory of Richard Paris}, SIGMA 20 (2024), 074.
\url{https://doi.org/10.3842/SIGMA.2024.074}

\bibitem{Hasanov-Ergashev-2023}
A. Hasanov, T.G. Ergashev, 
Infinite summation formulas for triple Lauricella hypergeometric functions, 
J. Math. Sci., New York 274 (2) (2023), 215--227.
\url{https://doi.org/10.1007/s10958-023-06590-z}

\bibitem{Hasanov-et al-2010}
A. Hasanov, M. Turaev, J. Choi, 
Decomposition formulas for the generalized hypergeometric ${}_{4}F_{3}$ fucnction, 
Honam Math. J. 32 (1) (2010), 1--16.  
\url{https://doi.org/10.5831/HMJ.2010.32.1.001}

\bibitem{Horn-1889}
J. Horn,
On the convergence of hypergeometric series of two and three variables,
Math. Ann. 34 (1889), 544--600.
\url{https://doi.org/10.1007/BF01443681} 

\bibitem{Karlsson-1986}
P.W. Karlsson, On intermediate Lauricella functions, J\~{n}\={a}n\={a}bha, 16 (1986), 211--222.

\bibitem{Komilova-et al-2024}
N.J. Komilova, A. Hasanov, T.G. Ergashev, 
Expansions of Kamp\'{e} de F\'{e}riet hypergeometric functions, 
J. Math., Mech. Comp. Sci. 122 (2024), 50--66. 
\url{https://doi.org/10.26577/JMMCS2024-122-02-b5}

\bibitem{Krantz-Book-2001}
S.G. Krantz, 
Function Theory of Several Complex Variables, 
AMS Chelsea Publishing, Providence, RI, 2001. Reprint of the 1992 edition.

\bibitem{Kristensson-Book-2010}
G. Kristensson,
Second Order Differential Equations: Special Functions and Their Classification,
Springer, New York, 2010.

\bibitem{Lauricella-1893}
G. Lauricella, 
Sulle funzioni ipergeometriche a pi\`{u} variabili, 
Palermo Rend. VII (1893), 111--158.
\url{https://doi.org/10.1007/BF03012437}

\bibitem{Luo-Raina-2021}
M.-J. Luo, R.K. Raina,
On certain results related to the hypergeometric function $F_K$, 
J. Math. Anal. Appl. 504 (2) (2021), Article ID 125439, 18 pp. 
\url{https://doi.org/10.1016/j.jmaa.2021.125439}

\bibitem{Luo et al-2022}
M.-J. Luo, M.-H. Xu, R.K. Raina, 
On certain integrals related to Saran's hypergeometric function $F_K$,
Fractal Fract. 6 (2022), 155.
\url{https://doi.org/10.3390/fractalfract6030155}

\bibitem{Olver et al-NISThandbook-2010}
F.W.J. Olver, D.W. Lozier, R.F. Boisvert, C.W. Clark (Eds.),
NIST Handbook of Mathematical Functions, Cambridge University Press, New York, 2010.

\bibitem{Pavlov-2023}
M. Pavlov, 
Global torus blocks in the necklace channel, 
Eur. Phys. J. C 83 (2023), Paper No. 1026, 12 pp.  
\url{https://doi.org/10.1140/epjc/s10052-023-12166-7}

\bibitem{Quan-2020}
J.-F. Quan, 
Study of multivariable hypergeometric series in conformal field theory (in Chinese), 
Bachelor's thesis, Donghua University, Shanghai, China, 2020.

\bibitem{Riemann-1857}
B. Riemann, 
Beitr\"{a}ge zur Theorie der durch die Gauss'sche Reihe $F(\alpha,\beta,\gamma,x)$ darstellbaren Functionen, 
Aus dem siebenten Band der Abhandlungen der Koniglichen Gesellschaft der Wissenschaften zu G\"{o}ttingen, 1857.

\bibitem{Rosenhaus-2019}
V. Rosenhaus,
Multipoint conformal blocks in the comb channel,
J. High Energy Phys. 2019 (2019), Paper No. 142, 35 pp.
\url{https://doi.org/10.1007/JHEP02(2019)142}

\bibitem{Sadykov-Tsikh-2014}
T.M. Sadykov, A.K. Tsikh, 
Hypergeometric and Algebraic Functions of Several Variables (Nauka, Moscow, 2014) [in Russian].

\bibitem{Saran-1954}
S. Saran, 
Hypergeometric functions of three variables,
Ganita 5 (2) (1954) 77--91.

\bibitem{Saran-1955a}
S. Saran,
Transformations of certain hypergeometric functions of three variables,
Acta Math. 93 (1955), 293--312.
\url{https://doi.org/10.1007/BF02392525}

\bibitem{Saran-1955b}
S. Saran,
Integrals associated with hypergeometric functions of three variables,
Proc. Natl. Inst. Sci. India, Part A 21 (1955) 83--90.

\bibitem{Saran-1956}
S. Saran,
The solutions of certain hypergeometric equations,
Proc. Natl. Inst. Sci. India, A Phys. Sci. 21 (1956) 404--408.

\bibitem{Semusheva-Tsikh-2009}
A.Y. Semusheva, A.K. Tsikh,
Domains of convergence of hypergeometric series of several complex variables,
J. Sib. Fed. Univ., Math. Phys. 2 (2) (2009), 221--229.

\bibitem{Srivastava-Karlsson-Book-1985}
H.M. Srivastava, P.W. Karlsson,
Multiple Gaussian Hypergeometric Series,
Halsted Press (Ellis Horwood Limited,
Chichester), Wiley, New York, Chichester, Brisbane and Toronto, 1985.

\bibitem{Srivastava-Manocha-Book-1984}
H.M. Srivastava, H.L. Manocha,
A Treatise on Generating Functions,
Ellis Horwood, Chichester, 1984.

\end{thebibliography}
\end{document}